\documentclass[12pt]{amsart}
\usepackage{
  fullpage,
  graphicx,amssymb,epic,eepic,
picins,
color,hyphenat,
  multicol, 
slashbox, 
}
\usepackage[all]{xy}
\usepackage{hyperref}

\theoremstyle{plain}
\newtheorem{theorem}{Theorem}[section]
\newtheorem{lemma}[theorem]{Lemma}
\newtheorem{proposition}[theorem]{Proposition}
\newtheorem{corollary}[theorem]{Corollary}

\theoremstyle{definition}
\newtheorem{remark}[theorem]{Remark}

\def\calO{{\mathcal O}}
\def\calI{{\mathcal I}}
\def\calA{{\mathcal A}}
\def\proj{{\operatorname{proj}\,}}
\def\bbGamma{{\raisebox{-0.11mm}{\begin{picture}(0,0)%
\includegraphics{figs/bbGamma.pstex}%
\end{picture}%
\setlength{\unitlength}{953sp}%
\begingroup\makeatletter\ifx\SetFigFont\undefined%
\gdef\SetFigFont#1#2#3#4#5{%
  \reset@font\fontsize{#1}{#2pt}%
  \fontfamily{#3}\fontseries{#4}\fontshape{#5}%
  \selectfont}%
\fi\endgroup%
\begin{picture}(464,606)(15729,-6933)
\end{picture}%
}}}

\def\arXiv#1{{\href{http://front.math.ucdavis.edu/#1}{arXiv:#1}}}

\catcode`\@=11
\long\def\@makecaption#1#2{%
    \vskip 10pt
    \setbox\@tempboxa\hbox{
      \small\sf{\bfcaptionfont #1. }\ignorespaces #2}%
    \ifdim \wd\@tempboxa >\captionwidth {%
        \rightskip=\@captionmargin\leftskip=\@captionmargin
        \unhbox\@tempboxa\par}%
      \else
        \hbox to\hsize{\hfil\box\@tempboxa\hfil}%
    \fi}
\font\bfcaptionfont=cmssbx10 scaled \magstephalf
\newdimen\@captionmargin\@captionmargin=2\parindent
\newdimen\captionwidth\captionwidth=\hsize
\catcode`\@=12

\begin{document}
\newdimen\captionwidth\captionwidth=\hsize

\title{Homomorphic Expansions for Knotted Trivalent Graphs}

\author{Dror Bar-Natan and Zsuzsanna Dancso}
\address{
  Department of Mathematics\\
  University of Toronto\\
  Toronto Ontario M5S 2E4\\
  Canada
}
\email{drorbn@math.toronto.edu, zsuzsi@math.toronto.edu}
\urladdr{http://www.math.toronto.edu/drorbn, http://www.math.toronto.edu/zsuzsi}

\date{First edition: Mar.~9, 2011. This edition: Jul.~21,~2012. {\em Electronic
versions:} \url{http://www.math.toronto.edu/~drorbn/papers/ktgs/} and
\arXiv{1103.1896}.}

\begin{abstract}
It had been known since old times \cite{MO, CL, Da} that there exists
a universal finite type invariant (``an expansion'') $Z^{old}$
for Knotted Trivalent Graphs (KTGs), and that it can be chosen to
intertwine between some of the standard operations on KTGs and their
chord-diagrammatic counterparts (so that relative to those operations,
it is ``homomorphic''). Yet perhaps the most important operation on KTGs
is the ``edge unzip'' operation, and while the behavior of $Z^{old}$
under edge unzip is well understood, it is not plainly homomorphic as some
``correction factors'' appear.

In this paper we present two (equivalent) ways of modifying $Z^{old}$
into a new expansion $Z$, defined on ``dotted Knotted Trivalent
Graphs'' (dKTGs), which is homomorphic with respect to a large set
of operations. The first is to replace ``edge unzips'' by ``tree
connected sums'', and the second involves somewhat restricting the
circumstances under which edge unzips are allowed. As we shall explain,
the newly defined class dKTG of knotted trivalent graphs retains all the
good qualities that KTGs have --- it remains firmly connected with the
Drinfel'd theory of associators and it is sufficiently rich to serve as a
foundation for an ``Algebraic Knot Theory''.  As a further application,
we present a simple proof of the good behavior of the LMO invariant
under the Kirby II (band-slide) move \cite{LMMO}.

\end{abstract}

\maketitle
\tableofcontents

\section{Introduction}\label{int}

Knot theory is not usually considered an algebraic subject, and one
reason for this is that knots are not equipped with a rich enough
algebraic structure. There are some operations defined on knots, 
most importantly connected sum and cabling, but even with these present,
the space of knots is far from finitely generated,
not to mention finitely presented. 

There is a way, however, to put knot theory in an algebraic context,
by considering the larger, richer space of Knotted Trivalent Graphs.
KTGs include knots and links, and are equipped with four standard
operations, called the orientation switch, edge delete, edge unzip,
and connected sum. With these operations, KTGs form a finitely
presented algebraic structure \cite{Th}. Furthermore, several topological
knot properties, including knot genus and the ribbon property are 
definable by simple formulas in the space of KTGs \cite{BN2}. Thus, invariants
which are well-behaved with respect to the algebraic structure on KTGs could
be used as algebraic tools to understand these knot properties.

A construction of an almost-perfect such invariant has long been known 
(\cite{MO}, and later \cite{CL} and even later \cite{Da}): the Kontsevich integral of knots can
be extended to a universal finite type invariant (or ``expansion'') $Z^{old}$
of KTGs, and the extension is very well-behaved with respect to three of
the four KTG operations. By ``very well-behaved'', we mean that it intertwines
those operations and their chord-diagrammatic counterparts, in other words,
it is ``homomorphic'' with respect to those operations. However, $Z^{old}$
fails to commute with the unzip operation, which plays a crucial role
in the finite generation of KTGs. Although the behavior of $Z^{old}$
with respect to unzip is well-understood, it is not homomorphic, and
it can be shown (we do so in the appendix
to this paper) that any expansion of KTGs will display an 
anomaly like the above: it can not commute with all four operations at the
same time.

The main goal of this paper is to fix the anomaly by proposing a
different definition of KTGs, which we will call ``dotted knotted trivalent
graphs'', or dKTGs on which a truly homomorphic expansion exists. We present
two (equivalent) constructions of this space. In one we replace the unzip, delete 
and connected sum operations by a more general set of operations called
``tree connected sums''. In the other, we restrict the set of edges which
we allow to be unzipped. We show that $Z^{old}$ can easily be modified to
produce a homomorphic expansion of dKTGs, and that dKTGs retain all the
good qualities of KTGs, namely, finite generation and a close connection
to Drinfel'd associators. We would like to emphasize that the value of this construction is 
not in extending $Z^{old}$ to dotted KTGs, but in providing a new way of
understanding and handling $Z^{old}$.

Finally, we show a simple (free of associators and local considerations) proof
of the theorem that the LMO invariant is well behaved with respect to the Kirby II
(band-slide) move.

\subsection{Acknowledgements}The authors would like to thank some
anonymous readers for many useful comments and for pointing out several
gaps in an earlier version of this paper.  We are grateful to the referees
at JKTR for their detailed and helpful suggestions.  The first author
was partially supported by NSERC grant RGPIN 262178 and the second author by
an NSERC PDF fellowship.

\section{Preliminaries}

The goal of this section is to introduce knotted trivalent graphs
and reintroduce the theory of finite type invariants
in a general algebraic context. We define
general ``algebraic structures'', ``projectivizations'' and ``expansions'',
with finite type invariants of 
knotted trivalent graphs (and the special case of knots and links) as an 
example.
 
\subsection{KTGs and $Z^{old}$}

\smallskip

A {\it trivalent graph} is a graph which has three edges meeting at each
vertex. We require that all edges be oriented
and that vertices be equipped with a cyclic orientation,
i.e. a cyclic ordering of the three half-edges meeting at the vertex. 
We allow multiple edges;
loops (i.e., edges that begin and end at the same vertex); and circles (i.e.
edges without a vertex).

\parpic[r]{\input{figs/thicken.pstex_t}}
Given a trivalent graph $\Gamma$, its {\it thickening} ${\raisebox{-0.11mm}{\begin{picture}(0,0)%
\includegraphics{figs/bbGamma.pstex}%
\end{picture}%
\setlength{\unitlength}{953sp}%
\begingroup\makeatletter\ifx\SetFigFont\undefined%
\gdef\SetFigFont#1#2#3#4#5{%
  \reset@font\fontsize{#1}{#2pt}%
  \fontfamily{#3}\fontseries{#4}\fontshape{#5}%
  \selectfont}%
\fi\endgroup%
\begin{picture}(464,606)(15729,-6933)
\end{picture}%
}}$ is
obtained from it by ``thickening'' vertices as shown on the right,
and gluing the resulting ``thick Y's'' in an orientation preserving
manner. Hence, the thickening is an oriented two-dimensional surface with
boundary. For it to be well-defined, we need the cyclic 
orientation at the vertices.

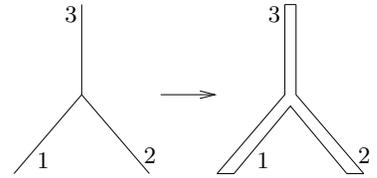
\parpic[r]{\begin{picture}(0,0)%
\includegraphics{figs/band.pstex}%
\end{picture}%
\setlength{\unitlength}{3355sp}%
\begingroup\makeatletter\ifx\SetFigFont\undefined%
\gdef\SetFigFont#1#2#3#4#5{%
  \reset@font\fontsize{#1}{#2pt}%
  \fontfamily{#3}\fontseries{#4}\fontshape{#5}%
  \selectfont}%
\fi\endgroup%
\begin{picture}(1866,1766)(1372,-2534)
\end{picture}%
}
A {\it Knotted Trivalent Graph (KTG)} is an embedding of a
thickened trivalent graph ${\raisebox{-0.11mm}{\begin{picture}(0,0)%
\includegraphics{figs/bbGamma.pstex}%
\end{picture}%
\setlength{\unitlength}{953sp}%
\begingroup\makeatletter\ifx\SetFigFont\undefined%
\gdef\SetFigFont#1#2#3#4#5{%
  \reset@font\fontsize{#1}{#2pt}%
  \fontfamily{#3}\fontseries{#4}\fontshape{#5}%
  \selectfont}%
\fi\endgroup%
\begin{picture}(464,606)(15729,-6933)
\end{picture}%
}}$ in $\mathbb{R}^3$, as shown.
The {\it skeleton} of a KTG $\gamma$ is the combinatorial 
object (trivalent graph $\Gamma$) behind it.
We consider KTGs up to isotopies that do not change the skeleton.
The thickening is equivalent to saying that the edges of the graph are
framed and the framings agree at vertices. In particular, framed knots and links
are knotted trivalent graphs. For a given trivalent graph $\Gamma$, let us denote
all the knottings of $\Gamma$ (i.e., all KTGs with skeleton $\Gamma$) by $\mathcal K(\Gamma)$.

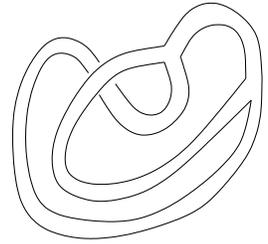
\parpic[r]{\begin{picture}(0,0)%
\includegraphics{figs/R1prime.pstex}%
\end{picture}%
\setlength{\unitlength}{4144sp}%
\begingroup\makeatletter\ifx\SetFigFont\undefined%
\gdef\SetFigFont#1#2#3#4#5{%
  \reset@font\fontsize{#1}{#2pt}%
  \fontfamily{#3}\fontseries{#4}\fontshape{#5}%
  \selectfont}%
\fi\endgroup%
\begin{picture}(954,1104)(3679,-3313)
\put(3916,-2716){\makebox(0,0)[lb]{\smash{{\SetFigFont{12}{14.4}{\rmdefault}{\mddefault}{\updefault}{\color[rgb]{0,0,0}$R1'$}%
}}}}
\end{picture}%
}
The reader may remember that framed knots are in one-to one correspondence
with knot diagrams modulo the usual Reidemeister moves $R2$ and $R3$, as
well as a weakened version of $R1$ that we call $R1'$, shown on the right.
($R1'$ does not change the framing but does change the winding number.)

\vspace{1cm}

Isotopy classes of KTG's are in one to one correspondence
with graph diagrams (projections onto a plane with only transverse double
points preserving the over- and under-strand information at the crossings),
modulo the Reidemeister moves $R1'$, $R2$, $R3$ and $R4$ (see for example \cite{MO}). 
We understand the framing 
corresponding to a graph diagram to 
be the blackboard framing. $R2$ and $R3$ are the same as in the knot
case. $R4$ involves moving a strand in front of or behind a vertex:
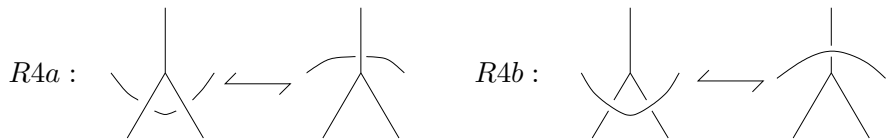
\begin{center}
\input figs/r4.pstex_t
\end{center}

There are four kinds of operations defined on KTG's:

\smallskip

Given a trivalent graph $\Gamma$, or a KTG $\gamma \in \mathcal K (\Gamma)$, 
and an edge $e$ of 
$\Gamma$, one can {\it switch the orientation} of $e$. We denote the
resulting graph by $S_e(\gamma)$. In other words, we have defined unary
operations $S_e: \mathcal K (\Gamma) \to \mathcal K (S_e(\Gamma))$.

\smallskip

We can also {\it delete} the edge
$e$, which means the two vertices at the ends of $e$ also cease to exist to
preserve the trivalence. To do this, it is required that the orientations of the
two edges connecting to $e$ at either end match. This operation
is denoted by $d_e: \mathcal K (\Gamma) \to \mathcal K (d_e(\Gamma))$.

\smallskip

{\it Unzipping} the edge $e$ (denoted by 
$u_e: \mathcal K (\Gamma) \to \mathcal K (u_e(\Gamma))$, see the figure below)
means replacing it by two edges that are ``very
close to each other''. The
two vertices at the ends of $e$ will disappear. 
This can be imagined as cutting the band of $e$ in 
half lengthwise. In the case of a trivalent graph $\Gamma$, we 
consider its thickening ${\raisebox{-0.11mm}{\begin{picture}(0,0)%
\includegraphics{figs/bbGamma.pstex}%
\end{picture}%
\setlength{\unitlength}{953sp}%
\begingroup\makeatletter\ifx\SetFigFont\undefined%
\gdef\SetFigFont#1#2#3#4#5{%
  \reset@font\fontsize{#1}{#2pt}%
  \fontfamily{#3}\fontseries{#4}\fontshape{#5}%
  \selectfont}%
\fi\endgroup%
\begin{picture}(464,606)(15729,-6933)
\end{picture}%
}}$ and similarly cut the edge $e$ in half
lengthwise.
Again, the orientations have to match,
i.e. the edges at the vertex where $e$ begins have to both be incoming,
while the edges at the vertex where $e$ ends must both be outgoing.

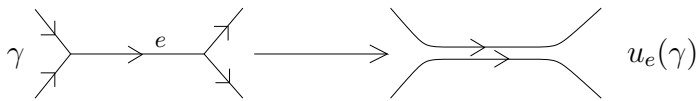
\begin{center}
\input figs/unzip.pstex_t
\end{center}

\smallskip

Given two graphs with selected edges $(\Gamma,e)$ and $(\Gamma', f)$,
the {\it connected sum} of these graphs along the two chosen edges,
denoted $\Gamma \#_{e,f} \Gamma'$, is
obtained by joining $e$ and $f$ by a new edge. For this to be well-defined,
we also need to specify the direction of the new edge,
the cyclic orientations at each new vertex, and 
in the case of KTGs, the framing on the new edge. 
To compress notation, let us declare that the new edge
be oriented from $\Gamma$ towards $\Gamma'$, have no twists,
and, using the blackboard framing, be attached to the right side
of $e$ and $f$, as shown:

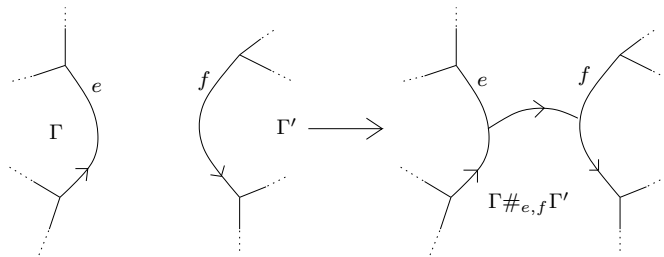
\begin{center}
\input figs/connsum.pstex_t
\end{center}

\begin{figure}
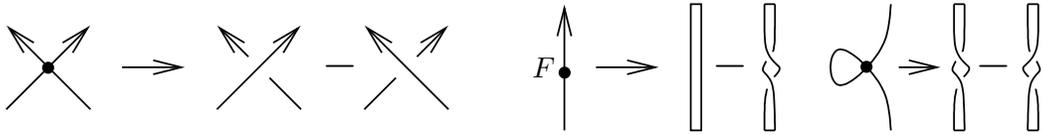

\[ \input figs/ftfilt.pstex_t \]
\caption{
  The resolution of double points and ``$F$-points''. (Note that
  ``$F$-points'' are needed because the resolution of a Reidemeister 1
  with double point involves two units of framing change rather than one.)
} \label{fig:ftfilt}
\end{figure}

The classical way of introducing finite type invariants for links,
extended straightforwardly to (framed) KTGs, is as follows: We allow formal 
linear
combinations of KTGs of the same skeleton, and
to filter the resulting vector space by resolutions
of singularities. An $n$-singular KTG is a trivalent graph 
immersed in $R^3$ with $n$ singular points: each may be a transverse double point
or merely a highlighted point on one of the strands marked by a letter ``$F$''. A resolution
of an $n$-singular KTG is obtained by replacing each
double point by the difference of an over-crossing and
an under-crossing, and each ``$F$-point'' by the difference between the 
graph obtained by its removal and by the graph obtained by replacing it 
with a 1-unit twist of framing (see Figure \ref{fig:ftfilt}).  
This produces a linear combination
of $2^n$ KTGs. The $n$-th piece of the finite type filtration
${\mathcal F}^n(\Gamma)$ is linearly generated by resolutions
of $n$-singular immersions of the skeleton $\Gamma$. Set 
$${\mathcal A(\Gamma)}=\bigoplus_{n=0}^{\infty}
{\mathcal F}^n(\Gamma)/{\mathcal F}^{n+1}(\Gamma),$$ the associated graded space
corresponding to the filtration. 
Finally, let ${\mathcal A}=\bigcup_{\Gamma}{\mathcal A}(\Gamma)$. 

\parpic[r]{\begin{picture}(0,0)%
\includegraphics{figs/ktgcd.pstex}%
\end{picture}%
\setlength{\unitlength}{3108sp}%
\begingroup\makeatletter\ifx\SetFigFont\undefined%
\gdef\SetFigFont#1#2#3#4#5{%
  \reset@font\fontsize{#1}{#2pt}%
  \fontfamily{#3}\fontseries{#4}\fontshape{#5}%
  \selectfont}%
\fi\endgroup%
\begin{picture}(1366,1370)(14798,-6504)
\end{picture}%
}
$\mathcal{A}(\Gamma)$ is best understood in terms of chord diagrams.
A \emph{chord diagram} of order $n$ on a skeleton graph $\Gamma$ is a combinatorial object
consisting of a pairing of $2n$ points on the edges of $\Gamma$,
up to orientation preserving homeomorphisms of the edges. Such
a structure is illustrated by drawing $n$ ``chords'' between the 
paired points, as seen in the figure on the right.
As in the finite type theory of links, a chord represents 
the difference of an over-crossing and an under-crossing
(i.e. a double point). This defines a map from chord diagrams
to ${\mathcal A}(\Gamma)$, which is well defined and surjective.

There are two classes of relations that are in the kernel,
the {\it Four Term relations ($4T$)}: 
\begin{center}
 \input figs/4t.pstex_t 
\end{center}
and the {\it Vertex Invariance relations ($VI$)},
(a.k.a. branching relation in \cite{MO}):
\begin{center}
\input figs/vi.pstex_t
\end{center}

In both pictures, there may be other chords in the parts of the graph not shown,
but they have to be the same throughout. In 4T, all skeleton parts (solid lines)
are oriented counterclockwise.
In $VI$, the sign $(-1)^{\rightarrow}$ is $-1$ if the edge the chord is ending
on is oriented to be outgoing from the vertex, and $+1$ if it is incoming.

Both relations arise from similar local isotopies of KTGs:

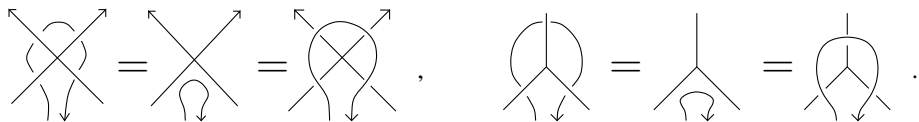
\begin{center}
\input{figs/trick.pstex_t}
\end{center}

Although it is easy to see that these relations are present,
showing that there are no more is difficult, and is best achieved 
by constructing 
a {\it universal finite type invariant}
$\mathbb Q KTG \to \mathcal A$ (we do not define universal finite type invariants here,
but will do so later in the general context). This was 
first done in \cite{MO} by extending the Kontsevich integral $Z$ of knots
\cite{Ko, CD, BN1}, 
building on results by T. Le, H. Murakami, J. Murakami and T. Ohtsuki
and using Drinfeld's theory of associators.  
In \cite{Da}, the same extension is constructed building on Kontsevich's original
definition. In this paper, we will denote this expansion by $Z^{old}$. 

Each operation of KTGs induces an operation on $\mathcal A$, as follows.

\smallskip

Given a graph $\Gamma$ and an edge $e$, the induced orientation switch
operation is a linear map $S_e:{\mathcal A}(\Gamma) \to {\mathcal A}(s_e(\Gamma))$
which multiplies a chord diagram $D$ by $(-1)^k$ where $k$ is the number of chords
in $D$ ending on $e$. 
Note that this generalizes the antipode map on Jacobi diagrams, which 
corresponds to the orientation reversal of knots (see \cite{Oh}, p.136).

\smallskip

The induced edge delete is a linear map
$d_e:{\mathcal A}(\Gamma) \to {\mathcal A}(d_e(\Gamma))$, defined as follows:
when the edge $e$ is deleted, all diagrams with a chord ending on $e$ are
mapped to zero, while those with no chords ending on $e$ are left unchanged,
except the edge $e$ is removed. 
Edge delete is the generalization
of the co-unit map of \cite{Oh} (p.136), and \cite{BN1}.

\smallskip

The induced unzip is a linear map
$u_e:{\mathcal A}(\Gamma) \to {\mathcal A}(u_e(\Gamma)).$
When $e$ is unzipped, each chord that ends on it is
replaced by a sum of two chords, one ending on each new edge (i.e., if
$k$ chords end on $e$, then $u_e$ sends this chord diagram
to a sum of $2^k$ chord diagrams). 

There is an operation on ${\mathcal A}(\circlearrowleft)$ corresponding to the cabling of knots:
references include \cite{BN1} (splitting map) and \cite{Oh} (co-multiplication). 
The graph unzip operation is
the graph analogy of cabling, so the corresponding map is analogous as well.

\smallskip

For graphs $\Gamma$ and $\Gamma'$, with edges $e$ and $e'$,
the induced connected sum
$\#_{e,e'}:{\mathcal A}(\Gamma) \times {\mathcal A}(\Gamma')
\to {\mathcal A}(\Gamma \#_{e,e'} \Gamma')$
acts in the obvious way, by performing the connected sum
operation on the skeletons and not changing the chords in any way. This is
well defined due to the $4T$ and $VI$ relations. 
(What needs to be proven
is that we can move a chord ending over the attaching point of the new edge;
this is done in the same spirit as the proof of Lemma 3.1 in \cite{BN1}, 
using ``hooks''; see also \cite{MO}, figure 4.)

\smallskip
\subsection{Algebraic structures and expansions}

An {\it algebraic structure} $\mathcal O$ is
some collection $(\mathcal O _\alpha)$ of sets of objects of different kinds,
where the subscript $\alpha$ denotes the {\it kind} of the objects in
$\mathcal O_ \alpha$, along with some collection of {\it operations} $\psi_\beta$,
where each $\psi_\beta$ is an arbitrary map with domain some product
$\mathcal O _{\alpha_1}\times\dots\times\mathcal O _{\alpha_k}$ of sets of objects,
and range a single set $\mathcal O _{\alpha_0}$ (so operations may be unary or
binary or multi-nary, but they always return a value of some fixed kind).
We also allow some named {\it constants} within some $\mathcal O _\alpha$'s
(or equivalently, allow some 0-nary operations).
The operations may or may not be subject to ``axioms'' --- an {\it axiom}
is an identity asserting that some composition of operations is equal
to some other composition of operations.

\parpic[r]{\begin{picture}(0,0)%
\includegraphics{figs/th.pstex}%
\end{picture}%
\setlength{\unitlength}{3789sp}%
\begingroup\makeatletter\ifx\SetFigFont\undefined%
\gdef\SetFigFont#1#2#3#4#5{%
  \reset@font\fontsize{#1}{#2pt}%
  \fontfamily{#3}\fontseries{#4}\fontshape{#5}%
  \selectfont}%
\fi\endgroup%
\begin{picture}(2424,924)(1189,-1573)
\end{picture}%
}
For example, KTGs form an algebraic structure that has a different kind 
of objects for each skeleton.
The sets of objects are the sets of knottings $\mathcal{K} (\Gamma)$ for each
skeleton graph $\Gamma$, and there are four operations, three unary and one binary.
KTGs are finitely generated\footnote{In the appropriate 
sense it is also finitely presented,
however we do not pursue this point here.}
(see \cite{Th}), by two
elements, the trivially embedded tetrahedron and the twisted tetrahedron,
shown on the right (note that these only differ in framing).

Any algebraic structure $\mathcal O$ has a {\it projectivization}. First extend
$\mathcal O$ to allow formal linear combinations of objects of the same kind
(extending the operations in a linear or multi-linear manner), then let
$\mathcal I$, the {\it augmentation ideal}, be the sub-structure made out of
all such combinations in which the sum of coefficients is $0$. Let
$\mathcal I^m$ be the set of all outputs of algebraic expressions (that is,
arbitrary compositions of the operations in $\mathcal O$) that have at least
$m$ inputs in $\mathcal I$ (and possibly, further inputs in $\mathcal O$), and
finally, set
$$
  \proj\calO:=\bigoplus_{m\geq 0} \calI^m/\calI^{m+1},
$$
the associated graded space with respect to the filtration by powers of the 
augmentation ideal.
Clearly, with the operation inherited from $\calO$, the projectivization
$\proj\calO$ is again an algebraic structure with the same multi-graph
of spaces and operations, but with new objects and with new operations
that may or may not satisfy the axioms satisfied by the operations of
$\calO$. The main new feature in $\operatorname{proj}\,\mathcal O$ is that it is a ``graded''
structure; we denote the degree $m$ piece $\mathcal I ^m/\mathcal I ^{m+1}$ of
$\operatorname{proj}\, \mathcal O$ by $\proj _m \mathcal O$. 

For example, in the case of KTGs, allow formal 
$\mathbb Q$-linear combinations of KTGs of the same skeleton and extend
the operations linearly.
The augmentation ideal $\mathcal I$ is generated by 
differences of knotted trivalent graphs of 
the same skeleton. $KTG$ is then filtered by 
powers of $\mathcal I$, and the projectivization 
$\mathcal A:=\proj KTG$ also
has a different kind of object for each skeleton $\Gamma$,
denoted $\mathcal A (\Gamma)$. 

\begin{theorem} \label{augm}
 The filtration by powers of the augmentation ideal $\calI$
coincides with the classical finite type filtration, and hence
the projectivization $\mathcal A$ defined above coincides with the space of
chord diagrams.
\end{theorem}

We defer the proof of this theorem to the appendix.

The finite type theory of knots and links is included in the 
finite type theory of KTGs
as a special case. On knots, there is no rich enough algebraic
structure for the finite type filtration to coincide with powers
of the augmentation ideal with respect to some operations. However,
knots and links form a subset of KTGs, and the restriction
of $\calI ^n$ to that subset reproduces the usual theory of 
finite type invariants of knots and links, and $Z^{old}$ restricts
to the Kontsevich integral. 

Given an algebraic structure $\mathcal O$
let $\operatorname{fil}\, \mathcal O$ denote the filtered structure of linear combinations of
objects in $\calO$ (respecting kinds), filtered by the powers $(\mathcal I^m)$
of the augmentation ideal $\mathcal I$. Recall also that any graded space
$G=\bigoplus_mG_m$ is automatically filtered, by $\left(\bigoplus_{n\geq
m}G_n\right)_{m=0}^\infty$.

An ``expansion'' $Z$ for $\mathcal O$
is a map $Z:\mathcal O \to \operatorname{proj}\, \mathcal O$ that preserves the kinds of objects
and whose linear extension (also called $Z$) to $\operatorname{fil}\, \mathcal O$ respects the
filtration of both sides, and for which $\left(\operatorname{gr}\, Z\right):
\left(\operatorname{gr}\, \operatorname{fil}\, \mathcal O = \operatorname{proj}\, \mathcal O \right) \to
\left(\operatorname{gr}\,\operatorname{proj}\,\mathcal O= \operatorname{proj}\, \mathcal O \right)$ is the identity map of
$\operatorname{proj}\, \mathcal O$.

In practical terms, this is equivalent to saying that $Z$ is a map
$\mathcal O \to \operatorname{proj}\, \mathcal O$ whose restriction to $\mathcal I^m$ vanishes in degrees
less than $m$ (in $\operatorname{proj}\, \mathcal O$) and whose degree $m$ piece is the
projection $\mathcal I ^m \to \mathcal I^m/ \mathcal I^{m+1}$.
Therefore, in the theory of KTGs, the word ``expansion'' means ``universal finite type invariant''.

A ``homomorphic expansion'' is an expansion which
also commutes with all the algebraic operations defined on the algebraic
structure $\mathcal O$.

\parpic[r]{\input{figs/glitch.pstex_t}}
For KTGs, it turns out (see \cite{MO, Da}) that $Z^{old}$ is almost homomorphic:
it intertwines the orientation switch, edge delete, and connected sum operations. 
However, $Z^{old}$ does not commute with edge unzip. The behavior with respect to unzip
is well-understood (showed in \cite{Da} using a result of \cite{MO}), and is described
by the formula shown in the figure on the right.
Here, $\nu$ denotes the Kontsevich integral of the unknot.
A formula for $\nu$ was conjectured in \cite{BGRT} and proven in \cite{BLT}.
The new chord combinations appearing on the right commute with all the old 
chord endings by 4T.
 A different way to phrase this formula is that
$Z^{old}$ intertwines the unzip operation 
$u_e: \mathcal{K} (\Gamma) \to \mathcal{K} (u_e(\Gamma))$ with a ``renormalized''
chord diagram operation $\tilde{u}_e: \mathcal{A}(\Gamma) \to \mathcal{A}(u_e(\Gamma))$,
$\tilde{u}_e=i_{\nu^{-1/2}}^2 \circ u_e \circ i_{\nu^{1/2}}$, where $i_{\nu^{1/2}}$
denotes the operation of placing a factor of $\nu^{1/2}$ on $e$, $u_e$ is the 
chord-diagram unzip operation induced by the topological unzip, and $i_{\nu^{-1/2}}^2$
places factors of $\nu^{-1/2}$ on each ``daughter edge''. So we have
$Z^{old}(u_e(\Gamma))=\tilde{u}_e Z^{old}(\Gamma)$.

This is an anomaly: if $Z^{old}$ was honestly homomorphic, 
there should be no new chords appearing, i.e., $Z$ should intertwine unzip and
its induced chord diagram operation.
Our main goal in this paper is to fix this.

\section{The homomorphic expansion}

Before we re-define $KTG$, let us note that doing so is really necessary:

\begin{theorem} \label{nonexistence}
 There is no homomorphic expansion $\mathbb{Q} KTG \to \mathcal{A}$, i.e.
an expansion cannot intertwine all four operations at once.
\end{theorem}

To keep an optimistic outlook in this paper, we defer the proof of this theorem to
the appendix.

\subsection{The space of dotted KTGs}

We define the algebraic structure $dKTG$ of 
{\it dotted Knotted Trivalent Graphs} as follows:

\parpic[r]{\begin{picture}(0,0)%
\includegraphics{figs/dktgex.pstex}%
\end{picture}%
\setlength{\unitlength}{4144sp}%
\begingroup\makeatletter\ifx\SetFigFont\undefined%
\gdef\SetFigFont#1#2#3#4#5{%
  \reset@font\fontsize{#1}{#2pt}%
  \fontfamily{#3}\fontseries{#4}\fontshape{#5}%
  \selectfont}%
\fi\endgroup%
\begin{picture}(609,591)(5275,-3250)
\end{picture}%
}
A {\it dotted trivalent graph} is a graph which may have
trivalent vertices, and two kinds of bivalent vertices, 
called dots, and anti-dots, the latter are
denoted by crosses. Dots and crosses are part of the skeleton
information. An example is shown on the right. 
Trivalent vertices are equipped 
with cyclical orientations and edges are oriented, as before.
Like before, a $dKTG$ $\Gamma$ has a well-defined
{\it thickening} ${\raisebox{-0.11mm}
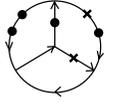
{\begin{picture}(0,0)%
\includegraphics{figs/bbGamma.pstex}%
\end{picture}%
\setlength{\unitlength}{953sp}%
\begingroup\makeatletter\ifx\SetFigFont\undefined%
\gdef\SetFigFont#1#2#3#4#5{%
  \reset@font\fontsize{#1}{#2pt}%
  \fontfamily{#3}\fontseries{#4}\fontshape{#5}%
  \selectfont}%
\fi\endgroup%
\begin{picture}(464,606)(15729,-6933)
\end{picture}%
}}$.

The algebraic structure $dKTG$ has a different kind of objects for
each dotted trivalent graph skeleton. The objects 
$\mathcal{K}(\Gamma)$ corresponding to skeleton
$\Gamma$ are embeddings of ${\raisebox{-0.11mm}{\begin{picture}(0,0)%
\includegraphics{figs/bbGamma.pstex}%
\end{picture}%
\setlength{\unitlength}{953sp}%
\begingroup\makeatletter\ifx\SetFigFont\undefined%
\gdef\SetFigFont#1#2#3#4#5{%
  \reset@font\fontsize{#1}{#2pt}%
  \fontfamily{#3}\fontseries{#4}\fontshape{#5}%
  \selectfont}%
\fi\endgroup%
\begin{picture}(464,606)(15729,-6933)
\end{picture}%
}}$ into 
$\mathbb{R}^3$, modulo ambient isotopy, or equivalently,
framed embeddings of $\Gamma$ where the framing agrees
with the cyclical orientations at trivalent vertices.

Obviously, $dKTG$s are also represented by $dKTG$
diagrams, with added Reidemeister moves to allow
the moving of bivalent vertices and anti-vertices
over or under an edge.

We define three kinds of operations on $dKTG$.

{\it Orientation reversal} reverses the orientation of an edge, as before. 

Given two $dKTG$s
$\gamma_1$ and $\gamma_2$ of skeletons $\Gamma_1$
and $\Gamma_2$, and with identical\footnote{By identical distinguished trees
we mean that the trees should be isomorphic
as graphs. Furthermore, the KTGs $\gamma_1$ and $\gamma_2$ are
each embedded inside unit cubes so that the distinguished trees
are the only parts positioned outside of the cube, to the right and
to the left respectively, forming mirror images of each other.
This can clearly be done by suitable isotopies of $\gamma_1$ and
$\gamma_2$.} distinguished 
trees $T_1$ and $T_2$, the
{\it tree connected sum} 
$\#_{T_1,T_2}: \mathcal K (\Gamma_1) \times \mathcal K (\Gamma_2)
\to \Gamma_1 \#_{T_1,T_2} \Gamma_2$
is obtained by deleting the two trees, and
joining corresponding ends by bivalent vertices, as shown.
The orientations of the new edges are inherited
from the leaves of the trees.
We leave it to the reader to check that this operation
is well-defined.

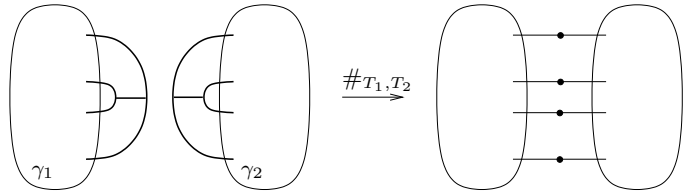
\begin{center}
 \input figs/tcs.pstex_t
\end{center}

We allow the distinguished trees to have dots and anti-dots on them, 
with the restriction that 
for each dot (resp. anti-dot) on $T_1$, $T_2$ is required to have 
an anti-dot (resp. dot) in the same position.

The {\it cancel} operation 
$c_{d,a}: \mathcal K (\Gamma) \to \mathcal K (c_{d,a} \Gamma)$
is defined when a $dKTG$ has a dot $d$ and an 
adjacent anti-dot $a$, in which case {\it cancel} 
deletes both (as in the figure below). 
This requires
the orientations of the three fused edges to agree.

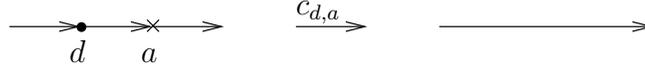
\begin{center}
 \input figs/cancel.pstex_t
\end{center}

\begin{lemma}\label{lem:fingen}
dKTGs with the above operations form a finitely generated algebraic structure.
\end{lemma}

\begin{proof}
We will show that orientation switch, unzip, delete and ``edge connected sum'' 
(that is, a connected sum operation 
followed by unzipping the connecting edge) are compositions of the new operations.
In the proof that $KTG$ is finitely generated (see \cite{Th}), a connected
sum is always followed by an unzip, so edge connected sum is sufficient for finite generation.
Furthermore, we need to show that it is possible to add dots and anti-dots
using tree connected sum, which then also allows one to delete any dots or anti-dots
using the cancel operation.

Orientation switch is an operation of $dKTG$s, so we have nothing to prove.

Unzip can be written as a tree connected sum the following way:

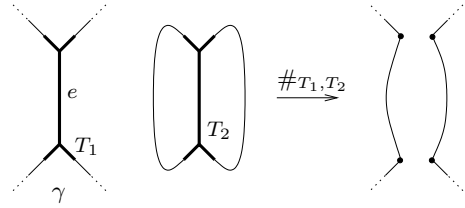
\begin{center}
 \input{figs/unziptcs.pstex_t}
\end{center}

\noindent
The graph on the right is almost $u_e(\gamma)$, except for the dots
which result from the tree connected sum. So to show that unzip can be 
written as a composition of the new operations, it is enough to show that it is
possible to ``get rid of'' dots. This is achieved by taking a tree
connected sum with a circle with three crosses and then canceling:

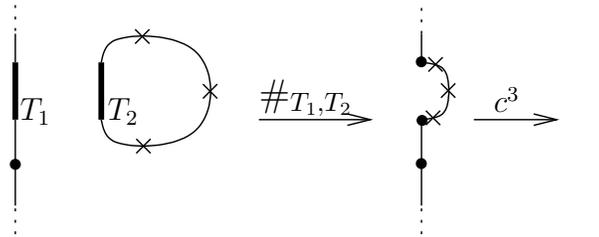
\begin{center}
 \input{figs/vtxkilling.pstex_t}
\end{center}

Edge delete and edge connected sum are done similarly, as illustrated by the figure below:

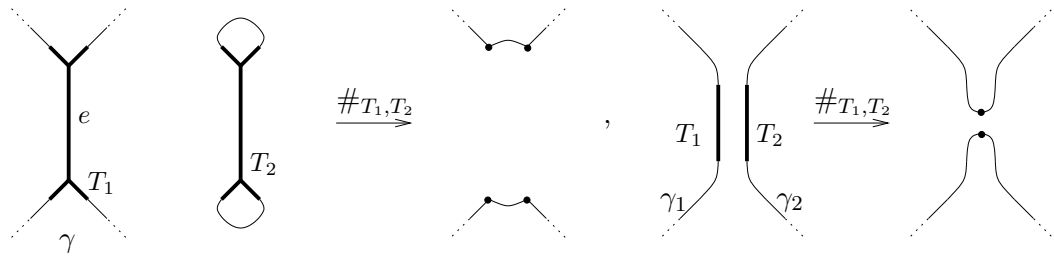
\begin{center}
 \input{figs/deletecs.pstex_t}
\end{center}
 
We have shown above that to add one anti-dot we need to take tree connected
sum with a circle with three anti-dots and a trivial tree (since the tree
connected sum produces two dots which then have to be canceled). Similarly,
to add one dot, we apply tree connected sum with a circle which has one 
anti-dot on it.

Using that $KTG$ is finitely generated by the two tetrahedrons,
we have now shown that $dKTG$ is finitely generated by the following four elements:

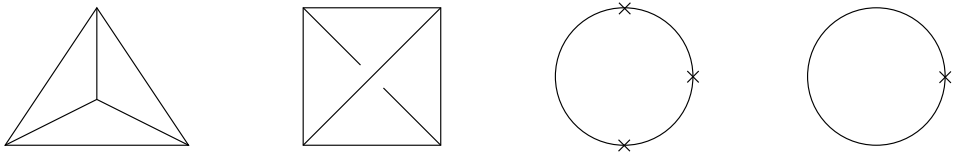
\begin{center}
 \begin{picture}(0,0)%
\includegraphics{figs/generators.pstex}%
\end{picture}%
\setlength{\unitlength}{3789sp}%
\begingroup\makeatletter\ifx\SetFigFont\undefined%
\gdef\SetFigFont#1#2#3#4#5{%
  \reset@font\fontsize{#1}{#2pt}%
  \fontfamily{#3}\fontseries{#4}\fontshape{#5}%
  \selectfont}%
\fi\endgroup%
\begin{picture}(6212,996)(1189,-1612)
\end{picture}%

\end{center}

\end{proof}

The reader might object that tree connected sum is infinitely many operations 
under one name, so it is not fair to claim that the structure is finitely generated. 
However, only two of these (the tree needed for unzip and delete, and the
trivial tree used for edge connected sum and vertex addition) are needed for finite 
generation. Later we will show a slightly different construction in which we only
use the operations that are essential, however, we felt that tree
connected sums are more natural and thus chose this version to be the starting point.

\subsection{The associated graded space and homomorphic expansion}

As before, the associated graded space has a kind of objects
for each skeleton dotted trivalent graph $\Gamma$, denoted $\mathcal A (\Gamma)$,
generated by chord diagrams on the skeleton $\Gamma$, and 
factored out by the
usual $4T$ and $VI$ relations, the latter of which now applies to dots 
and anti-dots as well, shown here for dots:

\begin{center}
\input{figs/newvi.pstex_t}
\end{center}

Orientation reversal acts the same way as it does for KTGs: if there are
$k$ chord endings on the edge that is being reversed, the diagram gets
multiplied by $(-1)^k$.

The tree connected sum operation acts on ${\mathcal A}^{dKTG}$ the following way:
if any chords end on the distinguished trees, we first use the VI relation to push 
them off the trees. Once the trees are free of chord endings, we join the skeletons
as above, creating bivalent vertices. Again, this operation is well-defined.

The cancel operation deletes a bivalent vertex and an anti-vertex on the same edge, 
without any change to chord endings.

\begin{theorem} \label{homexp}
There exists a homomorphic expansion $Z$ on the space of dKTGs, obtained 
from the $Z^{old}$ of \cite{MO} and \cite{Da} by placing a $\nu^{1/2}$ near each dot 
and $\nu^{-1/2}$ near each anti-dot.

\end{theorem}

\begin{proof}
First, note that $Z$ is well-defined: it does not matter which side of a dot (resp. anti-dot)
we place $\nu^{1/2}$ (resp. $\nu^{-1/2}$) on: if one edge is incoming, the other outgoing,
then these are equal by the $VI$ relation, otherwise they are equal by the $VI$ relation
and the fact that $S(\nu)=\nu$, where $S$ denotes the orientation switch operation.
(This is because the trivially framed oriented unknot is
isotopic to itself with the opposite orientation.)

Since $Z^{old}$ is an expansion of KTGs, it follows that $Z$ is an
expansion. For homomorphicity,
we must show that $Z$ commutes with the orientation switch, cancel and tree 
connected sum operations.

Orientation switch and cancel are easy. If an edge ends in two trivalent vertices,
then on that edge $Z$ coincides with $Z^{old}$ and hence commutes with switching 
the orientation. If one or both ends of the edge are bivalent, then $Z$ might
differ from $Z^{old}$ by a factor of $\nu$ (or two), but still commutes with $S$ by
the fact that $S(\nu)=\nu$. $Z$ commutes with cancel because the values of
the dot and anti-dot are inverses of each other (and local, therefore commute with all 
other chord endings and cancel each other out).

In terms of the ``old'' KTG operations (disregarding the dots for a moment, 
and ignoring edge orientation issues),
a tree connected sum can be realized by one ordinary connected sum followed by
a series of unzips:

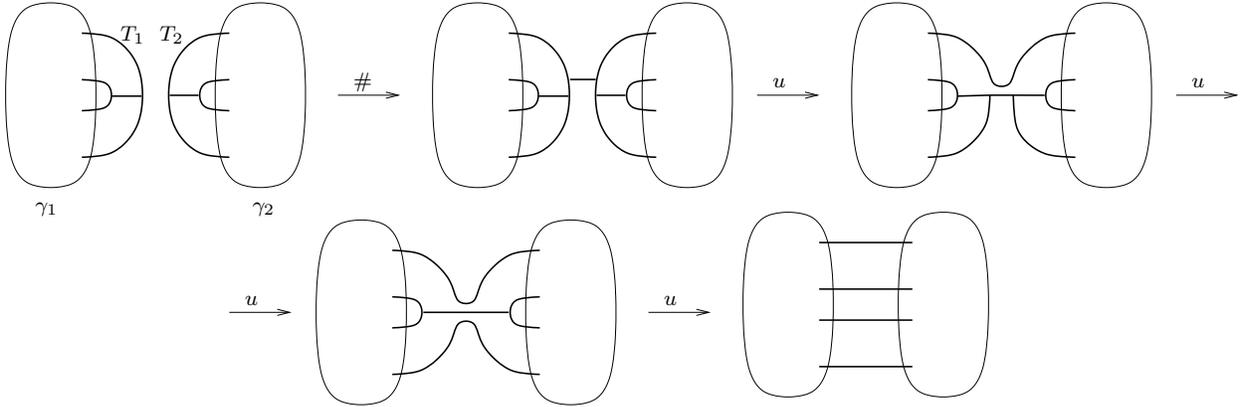
\begin{center}
 \input{figs/realizetcs.pstex_t}
\end{center}

We want to prove that $Z(\gamma_1 \#_{T_1,T_2} \gamma_2)=Z(\gamma_1)
\#_{T_1,T_2}Z(\gamma_2)$. To compute the left side, we trace $Z^{old}$
through the operations above. We assume that the trees have been cleared of
chord endings in the beginning using the $VI$ relation (and of course chords
that end outside the trees remain unchanged throughout). $Z^{old}$ commutes 
with connected sum, so in the first step, no chords appear on the trees.
In the second step, we unzip the bridge\footnote{We use the word ``bridge''
in the graph theoretic sense: an edge of the graph such that deleting this
edge would make its connected component break into two connected components.
In this context, we mean the new edge resulting from the connected sum.} 
connecting the two graphs. As mentioned
before, $Z^{old}$ intertwines unzip with 
$\tilde{u}=i_{\nu^{-1/2}}^2 \circ u \circ i_{\nu^{1/2}}$.
It is a simple fact of chord diagrams on KTGs that any chord diagram with a chord
ending on a bridge is $0$, thus the first operation $i_{\nu^{1/2}}$ is the identity
in this case. After the bridge is unzipped, $i_{\nu^{-1/2}}^2$ places $\nu^{-1/2}$
on the two resulting edges. The next time we apply $\tilde{u}$, $i_{\nu^{1/2}}$
cancels this out, the edge is unzipped, and then again $\nu^{-1/2}$ is placed
on the daughter edges, and so on, until there are no more edges to unzip. The 
operation $i_{\nu^{1/2}}$ will always cancel a $\nu^{-1/2}$ from a previous step. 
Therefore, at the end, the result is one factor of $\nu^{-1/2}$ on each of the
connecting edges.

We get $\gamma_1 \#_{T_1,T_2} \gamma_2$ by placing a dot on each of the connecting 
edges in the result of the above sequence of operations. $Z$ adds a factor of
$\nu^{1/2}$ at each dot, which cancels out each $\nu^{-1/2}$ that came from the
unzips. Thus, $Z(\gamma_1 \#_{T_1,T_2} \gamma_2)$ has no chords on the connecting
edges, which is exactly what we needed to prove.

Let us note that edge orientations can indeed be ignored: for the unzips
used above to be legitimate, a number of orientation switch operations
are needed, but since $S(\nu)=\nu$, the action of these on any chord diagram
that appears in the calculation above is trivial.

If the trees had dots and anti-dots to begin with, provided that for
every dot (resp.\ anti-dot) on $T_1$, $T_2$ had an anti-dot (resp.\
dot) in the same position, these will cancel each other out, and we have
already seen that $Z$ commutes with the cancel operation.
\end{proof}

\subsection{An equivalent construction} 
Let the space dKTG' have the same objects as dKTG, but we define the operations 
differently. We keep {\it orientation reversal} and {\it cancel} the same.
Instead of tree connected sums, we introduce the following three operations:

{\it Edge delete} is the same as in the space $KTG$, i.e., if orientations match, 
we can delete an edge connecting two trivalent vertices, and as a result, those
vertices disappear.

{\it Dotted unzip} allows unzips of an edge connecting two trivalent vertices
with one dot on it (technically, two edges), provided that orientations 
match at the trivalent vertices, as shown:

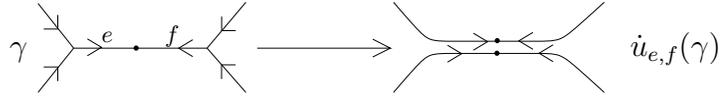
\begin{center}
 \input{figs/dotunzip.pstex_t}
\end{center}
  
{\it Dotted edge connected sum} is the same as edge connected sum, except 
dots appear where edges are ``fused'' (and there are no conditions on
edge orientations):

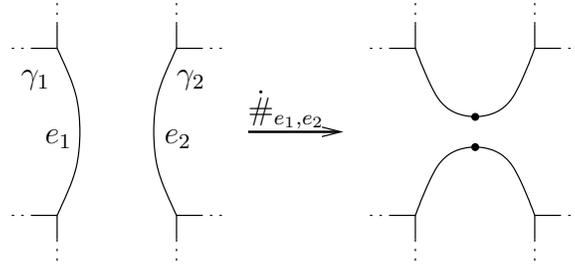
\begin{center}
 \input{figs/dotcs.pstex_t}
\end{center}

Alternatively, one can allow dotted connected sums (a connected sum where
a dot appears on the connecting edge). In this case, this construction is slightly
stronger than the previous one, as dotted connected sum cannot be written in
terms of the operations of the first construction. Dotted edge connected sum
is the composition of a dotted connected sum with a dotted unzip.

The associated graded space is as in the case of dKTG, and the induced operations
on it are the same for orientation reversal and cancel; are as in the case of
KTGs for delete and dotted unzip; and are as one would expect for dotted edge connected
sum (no new chords appear).

\begin{proposition}
The two constructions are equivalent in the sense that every dKTG operation can 
be written as a composition of dKTG' operations and every dKTG' operation is
a composition of dKTG operations, involving some constants (namely, a 
dumbbell graph and a trivially embedded theta graph, as shown in the proof).
\end{proposition}

\begin{proof}
For the first direction we only need to show that a tree connected sum can be
written as a composition of dKTG' operations. We have essentially done this
before, in the proof of Lemma \ref{lem:fingen}. The composition of operations required
is one dotted edge connected sum, followed by a succession of dotted unzips,
and orientation switches which we are ignoring for simplicity (as noted before,
they don't cause any trouble):

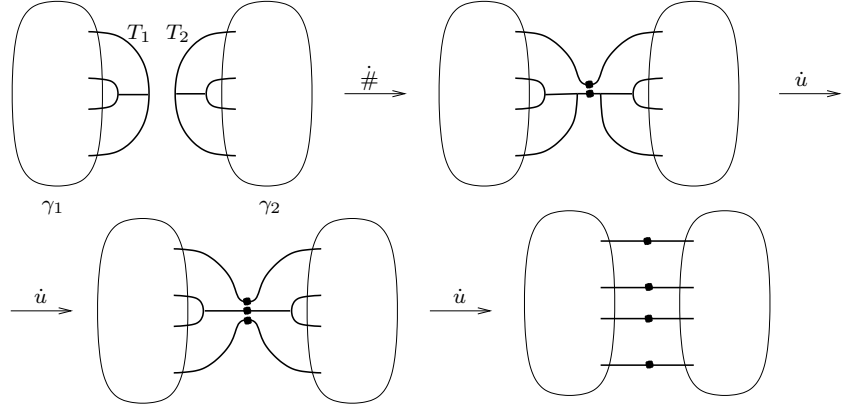
\begin{center}
 \input{figs/realizetcs2.pstex_t}
\end{center}

For the second direction, we need to write edge delete, dotted unzip and 
dotted edge connected sum as a composition of dKTG operations. Dotted
edge connected sum is a special case of tree connected sum (with trees consisting
of a single line segment).
Dotted unzip and edge delete are tree connected sums with given graphs and
given trees,
similar to the proof of Lemma \ref{lem:fingen}, shown below:
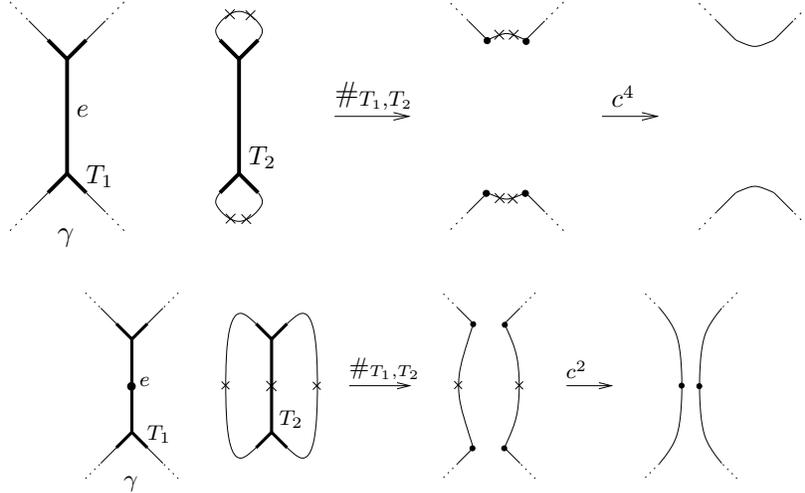
\begin{center}
\input{figs/gendelete.pstex_t}

\vspace{5mm}

\input{figs/genunzip.pstex_t}
\end{center}
\end{proof}

\begin{proposition}\label{hom}
 $Z$ is a homomorphic expansion of dKTG'.
\end{proposition}

\begin{proof}
It is obvious from the homomorphicity of $Z$ on dKTG that $Z$ 
commutes with orientation switch, cancel, and dotted edge connected sum.

Since dotted unzip and edge delete are unary operations, to show that
$Z$ commutes with them, we need to verify that the values of $Z$ 
on given graphs we used to produce these operations from tree connected sums 
are trivial, shown here for edge delete:
$$Z(d(\gamma))=Z(c^4(\gamma \# \gamma_0))=
c^4(Z(\gamma)\#Z(\gamma_0)).$$ 
Here, $\gamma_0$ denotes the ``dumbbell'' graph with four anti-dots,
shown above. If $Z(\gamma_0)=1$, and provided
that the edge to be deleted was cleared of chords previously, using the $VI$
relation, then the right side of the equation equals exactly $d(Z(\gamma))$.
Since $Z^{old}$ of the trivially embedded dumbbell (with no anti-dots) has a 
factor of $\nu$ on
each circle, $Z(\gamma_0)$ is indeed $1$, since two anti-dots add a factor
of $\nu^{-1}$ on each circle.

Unzip is done in an identical argument, where we use that the $Z$ value of
a trivially embedded theta-graph with one cross on each strand is 1. This
follows from the homomorphicity of $Z$: a tree connected sum of two of these 
graphs is again the graph itself, as shown below. If the $Z$ value of the graph 
is $\Psi$, we obtain that $\Psi^2=\Psi$, and since values of $Z$ are invertible
elements of $\mathcal A$, this implies that $\Psi=1$. 
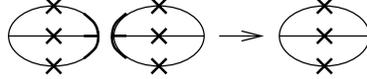
\begin{center}
 \begin{picture}(0,0)%
\includegraphics{figs/thetas.pstex}%
\end{picture}%
\setlength{\unitlength}{4144sp}%
\begingroup\makeatletter\ifx\SetFigFont\undefined%
\gdef\SetFigFont#1#2#3#4#5{%
  \reset@font\fontsize{#1}{#2pt}%
  \fontfamily{#3}\fontseries{#4}\fontshape{#5}%
  \selectfont}%
\fi\endgroup%
\begin{picture}(2184,494)(5029,-2558)
\end{picture}%

\end{center}
(As we will explain more precisely in Remark \ref{rem:spanningtree},
$\Psi$ here lives in $\mathcal A(\uparrow_2)$: chord diagrams on two
vertical strands, which have a natural algebra structure by concatenation
of the strands.)
\end{proof}

\section{The relationship with Drinfel'd associators}\label{sec:assoc}
Associators are useful and intricate gadgets that were first
introduced and studied by Drinfel'd \cite{Dr1, Dr2}.
The theory was later put in the context of parenthesized 
(a.k.a. non-associative) braids by \cite{LM}, \cite{BN3} and \cite{BN4}.
Here we present a construction of an associator as the value
of $Z$ on a dotted KTG.\footnote{This, along with the fact \cite{MO}
that $Z^{old}$ can be constructed from associators, is the content
of the often stated assertion that ``associators are equivalent
to homomorphic expansions of knotted trivalent graphs''.}

Let us first remind the reader of the definitions. 
Let $\calA(\uparrow_n)$ denote the space of chord diagrams on
$n$ upward oriented vertical strands subject to the $4T$ relations. Note that $\calA(\uparrow_n)$
has a natural algebra structure by stacking (concatenation) of the strands.
An associator is an element $\Phi \in \calA(\uparrow_3)$, which satisfies
three major equations, called the pentagon and the two hexagon equations,
as well as two minor conditions.

The pentagon is an equation in $\calA(\uparrow_4)$. Before we write it,
let us define some necessary maps $\Delta_i: \calA(\uparrow_n) \to \calA(\uparrow_{n+1})$.
$\Delta_i$, for $i=1,2,..,n$, is the doubling (unzip) of the $i$-th strand, which
acts on chord diagrams the same way unzip does. $\Delta_0$ adds an empty strand on the
left, leaving chord diagrams unchanged. Similarly, $\Delta_{n+1}$ adds a strand on the right.
Multiplication in $\calA(\uparrow_n)$ is defined by stacking chord diagrams
on top of each other.
In this notation, we can write the pentagon equation as follows
(products are read left to right and bottom to top):
$$\Delta_4(\Phi) \cdot \Delta_2(\Phi) \cdot \Delta_0(\Phi)=
\Delta_1(\Phi) \cdot \Delta_3(\Phi).$$

\parpic[r]{\input{figs/permute.pstex_t}}
The hexagons are two equations in $\calA(\uparrow_3)$, involving $\Phi \in \calA(\uparrow_3)$ and an extra
$R \in \calA(\uparrow_2)$. The permutation group $S_n$ acts on 
$\calA(\uparrow_n)$ by permuting the strands. We denote this action by superscripts,
as illustrated by the example on the right. The two hexagon equations are as follows:

$$\Phi \cdot \Delta_2(R) \cdot \Phi^{231} = \Delta_3(R) \cdot \Phi^{213} \cdot (\Delta_0(R))^{213}$$
$$\Phi \cdot \Delta_2((R^{21})^{-1}) \cdot \Phi^{231} = 
\Delta_3((R^{21})^{-1}) \cdot \Phi^{213} \cdot \big(\Delta_0((R^{21})^{-1})\big)^{213}$$

In addition to the pentagon and hexagon equations, associators are required to satisfy two minor conditions, 
non-degeneracy and horizontal unitarity (or mirror skew-symmetry), as below:
\begin{enumerate}
 \item $\Phi$ is {\it non-degenerate}: $d_i(\Phi)=1$ for $i=1,2,3$, where $d_i$ denotes the deletion
of the $i-th$ strand, which acts the same way as it does for KTGs, and maps 
$\mathcal{A}(\uparrow_3)$ to $\calA(\uparrow_2)$.
 \item The mirror image of $\Phi$ is $\Phi^{-1}$: $\Phi^{321}=\Phi^{-1}$.
\end{enumerate}

There are several other useful symmetry properties that the associator constructed here satisfies. Some of
these are required in constructions which involve associators as ingredients, for example, 
the construction of a universal finite type invariant of KTGs in \cite{MO} uses the first two properties
mentioned below. Stating some of these is notationally awkward, so let us introduce some useful 
alternative notation before defining the properties.

For any map between free groups $\beta: F_m \to F_n$, there exists a ``pullback'' $\beta^*: \calA_n \to \calA_m$,
defined the following way. $\beta$ defines a map of a bouquet of $m$ circles to
a bouquet of $n$ circles, up to homotopy, which induces $\beta$ on the fundamental groups. Removing the
joining point of the bouquets, one obtains a ``partial covering'' map from $m$ oriented lines to $n$ oriented lines.
Given a chord diagram $c \in \calA_n$, $\beta^*(c)$ is the sum of all possible lifts to the covering $m$ lines,
with negative signs wherever the orientations are opposite, and setting the image to be $0$ when a chord
ends on a strand not covered. Let us illustrate this on an example.

\parpic[r]{\input{figs/ex1.pstex_t}}
Consider the map $\beta=(x_1, x_1x_2^{-1}, x_2^{-1}): F_3 \to F_2$, where $x_1$ and $x_2$
are the free generators of $F_2$, and the notation means that $\beta$ sends the first free generator 
of $F_3$ to $x_1$, the second to $x_1x_2^{-1}$, and the third to $x_2^{-1}$. In the covering sense, this
can be illustrated by the left side of the figure on the right, while the right side shows the image
$\beta^*(c)$ of a chord diagram $c \in \calA_2$ in familiar, but awkward notation.

Properties $(1)$ and $(2)$ are easily rephrased in terms of the new notation:
\begin{enumerate}
\item $d_1(\Phi)=1$ can be written as $\beta_1^*(\Phi)=1$ where $\beta_1=(x_2, x_3):F_2 \to F_3$. 
Similarly for $d_2$ and $d_3$. 
\item $\Phi^{321}=\Phi^{-1}$ means $\beta_2^*(\Phi)=\Phi^{-1}$, where 
$\beta_2=(x_3,x_2,x_1):F_3 \to F_3$.
\end{enumerate}

Now let us state the further symmetry properties promised above.

\begin{enumerate}
 \setcounter{enumi}{2}
 \item \parpic[r]{\input{figs/horproperty.pstex_t}}
\noindent
An associator is called {\it horizontal} if it lives in the sub-algebra \linebreak $\calA^{hor}(\uparrow_3)$,
consisting of chord diagrams made only of horizontal chords. The associator we construct here is not a horizontal
chord associator, however, it possesses the following
crucial property of horizontal chord associators: $\beta_3^*(\Phi)=1$, where 
\linebreak $\beta_3=(x_2x_1^{-1},x_3x_1^{-1}): F_2 \to F_3$. 
This, in fact, is the only consequence of horizontality needed in \cite{MO}. The figure on the right shows this 
property in the other notation.
 
 \item \label{def:evenness} $\Phi$ is {\it even} or {\it vertically unitary}, meaning $\Phi$ ``upside down'' is 
$\Phi^{-1}$: $S_1S_2S_3(\Phi)=\Phi^{-1}$, equivalently $\beta_4^*(\Phi)=\Phi^{-1}$, where 
$\beta_4=(x_1^{-1},x_2^{-1},x_3^{-1}):F_3 \to F_3$.
 
\item \parpic[r]{\input{figs/rotsymm.pstex_t}}
\noindent
We will prove that the associator we construct has an additional symmetry, as shown
in the figure below, which we call {\it rotational symmetry}: 
$\beta_5^*(\Phi)=\Phi$, where \linebreak $\beta_5=(x_3^{-1},x_3^{-1}x_1,x_2^{-1}x_1)$. This
arises from the order $3$ rotational symmetry of the tetrahedron, and the reader can check that indeed
$\beta_5^3=1$. This property is rather unsightly when written in the other notation, nevertheless it is
shown on the right.
\end{enumerate}

\parpic[r]{\raisebox{-30mm}{\input{figs/phiiso.pstex_t}}} \picskip{6}
\begin{remark}\label{rem:spanningtree}
Note that for any chord diagram on a dKTG skeleton, one can pick any spanning tree and use the $VI$ 
relation to ``sweep it free of chords''. (In a slight abuse of notation, by a vertex we shall mean a
trivalent vertex, and by an edge, and edge connecting two trivalent vertices, which may have dots and
crosses on it, so it may really be a path.) This ``sweeping trick'' induces a (well-defined) isomorphism
from chord diagrams on a dKTG with a specified spanning tree to some $\calA(\uparrow_n)$. For example, there 
is an isomorphism from chord diagrams on a trivially embedded tetrahedron to $\calA(\uparrow_3)$, as shown 
in the figure on the right.
\end{remark}

We can now state and prove that we have constructed an associator:
\begin{theorem}\label{assoc}
The following $\Phi$ and $R$ satisfy the pentagon and hexagon equations, as well as properties $(1)$ -- $(5)$, 
and therefore $\Phi$ is a (nice) associator:
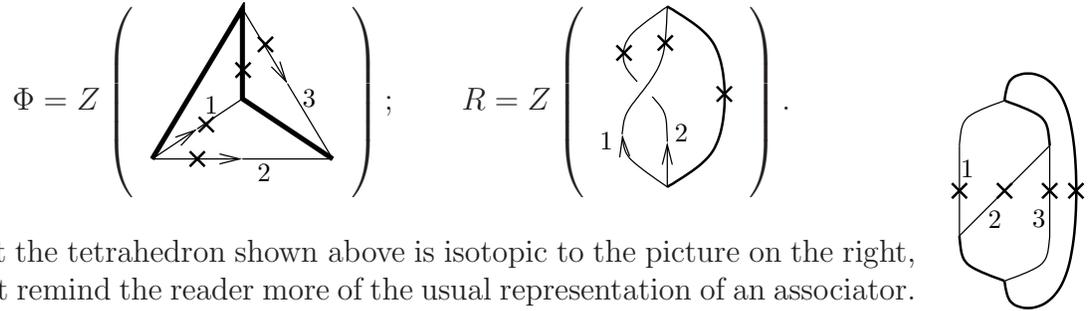
\[
  \Phi=Z\left(\begin{array}{c}\input{figs/phi.pstex_t}\end{array}\right);
  \qquad
  R=Z\left(\begin{array}{c}\input{figs/r.pstex_t}\end{array}\right).
  \qquad\qquad\null
\]
\end{theorem}

\parpic(0pt,0pt)(0cm,1cm)[r]{\input{figs/otherphi.pstex_t}}\picskip{3}
Note that the tetrahedron shown above is isotopic to the picture on the right, which might remind the reader
more of the usual representation of an associator.

\begin{remark}
The reader may wonder if this is a new construction for an associator
$\Phi$.  The answer is that it somewhat depends on how we pick the
invariant $Z^{old}$ at the start of our discussion. If we take it from
\cite{MO}, we must note that \cite{MO}'s construction in itself starts
with an associator, which we will call here $\Phi_0$. One may show that if
$\Phi_0$ is even (or vertically unitary, see (\ref{def:evenness}) above),
then $\Phi=\Phi_0$ and we got nothing new. But even if $\Phi_0$ is not
vertically unitary, our $\Phi$ here is vertically unitary (and satisfies
all other properties (1)-(5) above), so in that case $\Phi\neq\Phi_0$. Yet
note that in the process of constructing $\Phi$ from $\Phi_0$ we multiply
by several normalization factors and hence we cannot guarantee that $\Phi$
will consist only of horizontal chords. Finally, if we take $Z^{old}$ from
\cite{Da}, then no associator is used explicitly within its construction,
yet the normalization procedures employed in \cite{Da} are equivalent to
an implicit construction of the standard ``KZ'' associator $\Phi_{KZ}$
followed by the \cite{MO} procedure with seed $\Phi_0=\Phi_{KZ}$. Since
$\Phi_{KZ}$ is not vertically unitary, in that case $\Phi\neq\Phi_{KZ}$,
$\Phi$ satisfies (1)-(5) above, yet $\Phi$ is not necessarily made
entirely of horizontal chords.
\end{remark}

\parpic[r]{\begin{picture}(0,0)%
\includegraphics{figs/vconn.pstex}%
\end{picture}%
\setlength{\unitlength}{4144sp}%
\begingroup\makeatletter\ifx\SetFigFont\undefined%
\gdef\SetFigFont#1#2#3#4#5{%
  \reset@font\fontsize{#1}{#2pt}%
  \fontfamily{#3}\fontseries{#4}\fontshape{#5}%
  \selectfont}%
\fi\endgroup%
\begin{picture}(1752,834)(5119,-2683)
\put(5986,-2219){\makebox(0,0)[b]{\smash{{\SetFigFont{10}{12.0}{\rmdefault}{\mddefault}{\updefault}{\color[rgb]{0,0,0}$\#$}%
}}}}
\end{picture}%
}
{\it Proof.}
We first prove that $\Phi$ satisfies the pentagon equation. In the proof, we will use the
``vertex connected sum'' operation shown in the
figure on the right. This can be thought of either as a tree connected sum in the first model, or the composition
of a dotted connected sum composed with dotted unzips in the second model.

Now let us consider the following sequence of operations. In each step, we are thinking of $Z$ of the
pictured graph, which commutes with all the operations. To save space, we will not write out the $Z$'s.
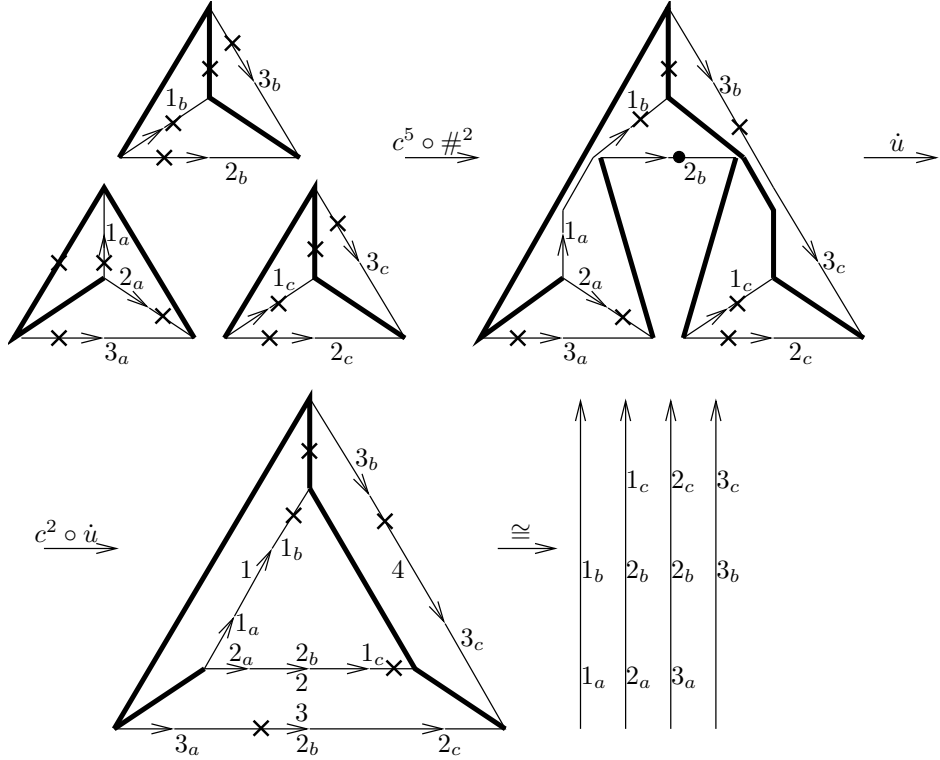
\begin{center}
 \input{figs/lspent.pstex_t}
\end{center}
Since $Z$ is homomorphic, the result of this sequence of operations is
$\Delta_4(\Phi) \cdot \Delta_2(\Phi) \cdot \Delta_0(\Phi),$
the left side of the pentagon equation.

For the right side, we perform a vertex connected sum of two tetrahedra:
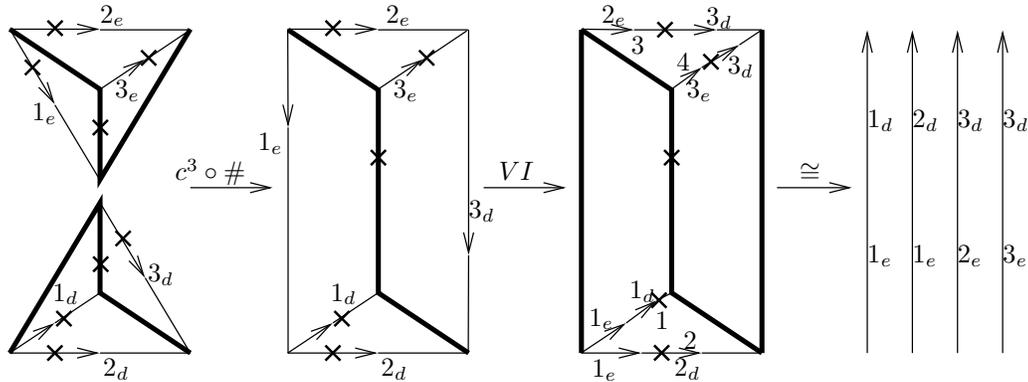
\begin{center}
 \input{figs/rspent.pstex_t}
\end{center}
The result can be written as $\Delta_1(\Phi) \cdot \Delta_3(\Phi)$, the right side
of the pentagon equation.

Since the two resulting dKTGs are isotopic (trivially embedded triangular prisms
with crosses in the same positions), the
two results have to be equal, and therefore $\Phi$ satisfies the pentagon equation.

Morally, the hexagon equation amounts to adding a twist to one of the tetrahedra in
the right side of the pentagon, on the middle crossed edge, which produces a triangular
prism with a twist on the middle vertical edge. Unzipping this edge then gives a new
twisted tetrahedron. 

More precisely, we carry this out in a similar fashion to the proof of the pentagon. 
For the left side 
of the first hexagon equation we take vertex connected sum
of two tetrahedra with a twisted theta graph:
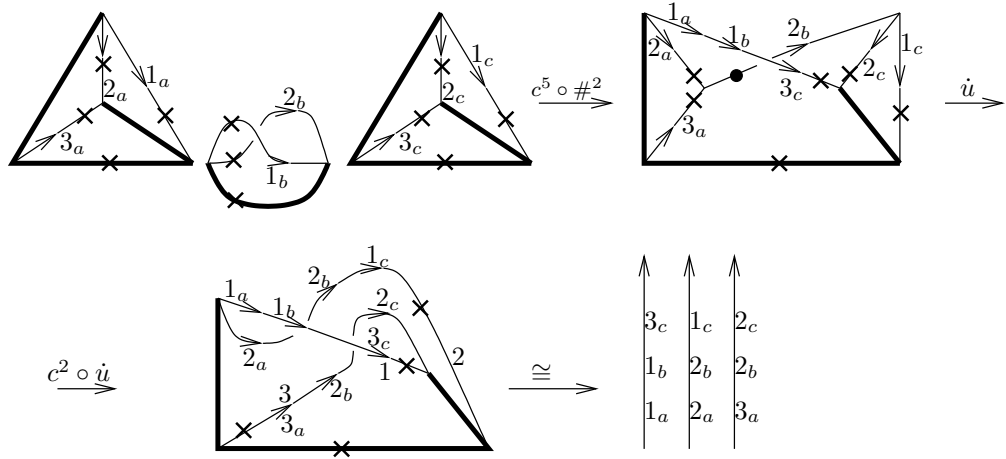
\begin{center}
 \input{figs/lshex.pstex_t}
\end{center}
The result is $\Phi \cdot \Delta_2(R) \cdot \Phi^{231}$, the left side of the first hexagon equation.

For the right side, we need to add an extra strand to the twisted theta graphs. This
can be done without changing the value of $Z$: we have mentioned before
that the $Z$- value of a trivially embedded theta graph with one cross on each strand 
is 1. Adding an extra strand can be realized by taking an edge connected sum with 
such a theta graph:
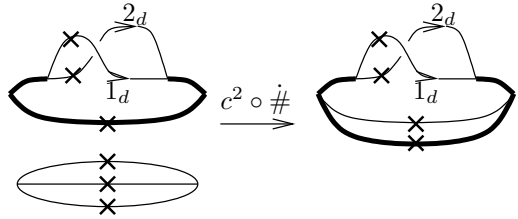
\begin{center}
 \input{figs/addstrand.pstex_t}
\end{center}
 
On to the right side of the hexagon equation, we now connect a tetrahedron with two twisted theta
graphs with an added strand, and unzip. The reader can check that moving the extra strand from 
one side to the other can be done by an isotopy.
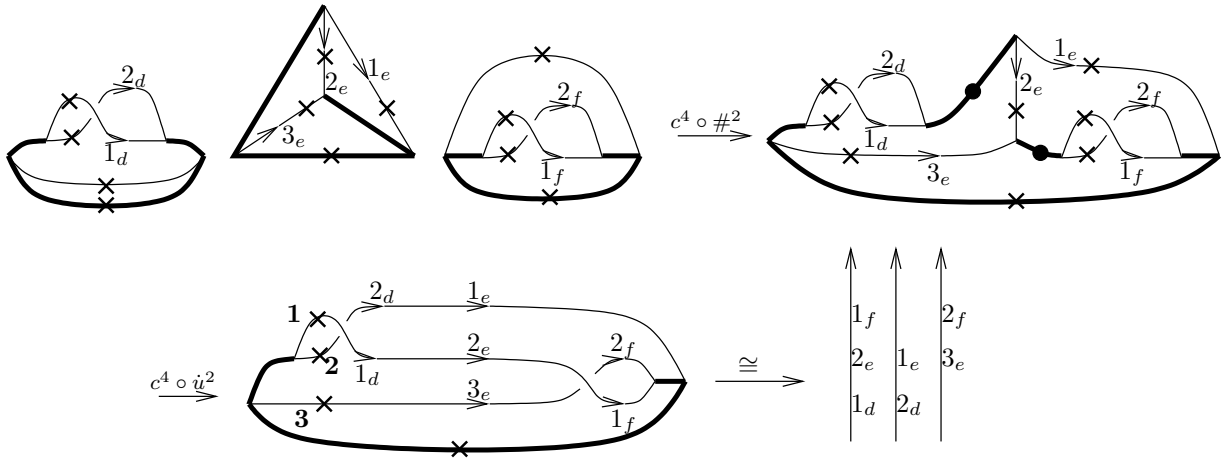
\begin{center}
 \input{figs/rshex.pstex_t}
\end{center}
The result reads $\Delta_3(R) \cdot \Phi^{213} \cdot (\Delta_0(R))^{213}$, the right hand
side of the hexagon. Is is apparent that the resulting twisted tetrahedron graphs of the
left and right sides are isotopic, proving that $\Phi$ and $R$ satisfy the first hexagon
equation.

For the second hexagon, we first show that
$$(R^{21})^{-1}=\raisebox{-7mm}{\input figs/r21inverse.pstex_t},$$ by
taking a vertex connected sum with $R$:
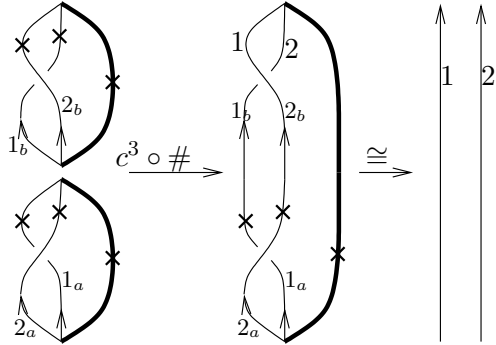
\begin{center}
 \input{figs/rproof.pstex_t}
\end{center}
Note that the resulting graph is isotopic to a trivially embedded theta-graph with one cross on
each edge. We have seen at the end of the proof of Proposition \ref{hom} that $Z$ evaluated
on this graph is trivial, which proves the claim. For a proof of the second hexagon equation, one
substitutes the picture for $(R^{21})^{-1}$ in place of $R$ everywhere.

\begin{remark} It is easy to check from the definition of $Z$ (for
details see \cite[Chapter 5.1]{Da2}), that $R^{21}=R$. If one requires
this additional property of $R$, then the strand swap of the second
hexagon equation may be omitted.  In fact $R$, when thought of as a chord
diagram on two strands, has only horizontal chords.  More specifically,
$R=e^{\frac{t_{12}}{2}}$, where $t_{12}$ is a chord connecting the
two strands.  See also Aside \ref{nonoriented}.
\end{remark}

Now on to proving the additional symmetry properties:

\begin{enumerate}
 \item To show non-degeneracy, observe that by homomorphicity 
$$d_i(\Phi)=Z\Bigg(d_i\Bigg( \raisebox{-9mm}{\input{figs/phi.pstex_t}} \Bigg) \Bigg).$$
Deleting any of the numbered edges from the tetrahedron yields a theta-graph
with one cross on each edge, and, as mentioned before, the value of $Z$ on this
graph is indeed $1$.

\item 
Horizontal unitarity follows from the fact that $\Phi$ satisfies the hexagon equations and
that $R^{21}=R$, as proven in \cite{BN3} Proposition 3.7, or \cite{Dr1} Proposition 3.5.

Alternatively, observe that $\Phi^{321}$ is the horizontal mirror image of $\Phi^{-1}$,
and $Z$ is invariant under taking mirror images for dKTGs with no crossings. This follows 
from the definition of $Z^{old}$, which is based on the original construction of the
Kontsevich integral: for details see \cite{Da2}, Chapter 5.1.

\item The horizontal-like property has the following pictorial proof:
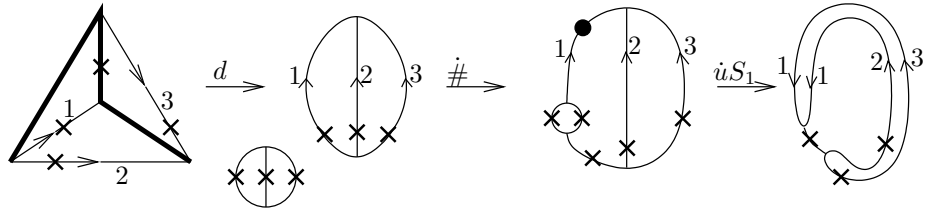
\begin{center}
 \input{figs/horiz.pstex_t}
\end{center}
In the first step, we are deleting the un-numbered edge which carries a cross.
The second step doesn't change the $Z$-value of the graph (as an element of 
$\calA(\uparrow_3)$), since we are taking a connected sum with a dKTG with
a trivial $Z$-value. The rest follows from homomorphicity. The $Z$-value of the graph on the right is,
on one hand, the same as required in the statement, as an element of 
$\calA(\uparrow_2)$, on the other hand, it is the theta-graph with three crosses,
which has trivial value.

\item Evenness follows from unitarity: the graphs representing $S_1S_2S_3(\Phi)$ and $\Phi^{321}$
are isotopic.

\item
And finally, rotational symmetry is due to the rotational symmetry of the tetrahedron. The right hand
side of the figure below is the same as the formula in the statement.
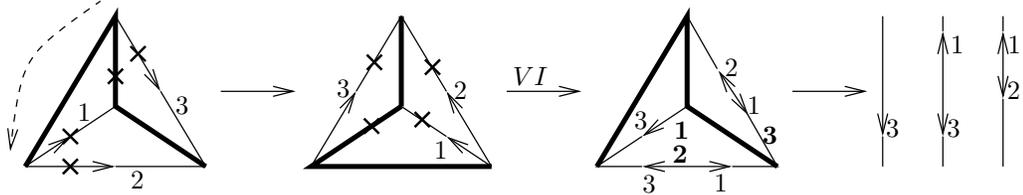
\begin{center}
 \input{figs/rotation.pstex_t}
\end{center}
\end{enumerate}
\qed

\subsection{Aside}\label{nonoriented} 
As a small detour, we note that $R^{21}=R$ could also be deduced from the properties of a homomorphic expansion,
if defined slightly more generally. In this paper we work with KTGs that are oriented surfaces embedded
in $R^3$. However, one could allow non-oriented surfaces, with a slight modification of the methods presented
here. In this case, the Moebius band replaces the twisted tetrahedron as a generator of the space of dKTGs.

\parpic[r]{\begin{picture}(0,0)%
\includegraphics{figs/nonorientedex.pstex}%
\end{picture}%
\setlength{\unitlength}{4144sp}%
\begingroup\makeatletter\ifx\SetFigFont\undefined%
\gdef\SetFigFont#1#2#3#4#5{%
  \reset@font\fontsize{#1}{#2pt}%
  \fontfamily{#3}\fontseries{#4}\fontshape{#5}%
  \selectfont}%
\fi\endgroup%
\begin{picture}(1037,894)(5434,-3268)
\end{picture}%
}
The notation needs to be modified as blackboard framed graph diagrams are not enough in this case:
we need an additional feature, the ``half twist'', which comes in left and right varieties, and the vertex
orientations need to be specified at each vertex (see the example on the right). As for Reidemeister moves,
the half twists can move freely behind and in front of crossings and across dots and crosses. 
There are a few additional relations: one can move a
twist across a vertex, two twists equal a kink, and two opposite twists cancel, shown below for one direction:
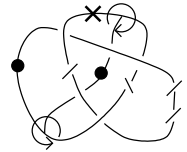
\begin{center}
 \begin{picture}(0,0)%
\includegraphics{figs/nonorientedrels.pstex}%
\end{picture}%
\setlength{\unitlength}{4144sp}%
\begingroup\makeatletter\ifx\SetFigFont\undefined%
\gdef\SetFigFont#1#2#3#4#5{%
  \reset@font\fontsize{#1}{#2pt}%
  \fontfamily{#3}\fontseries{#4}\fontshape{#5}%
  \selectfont}%
\fi\endgroup%
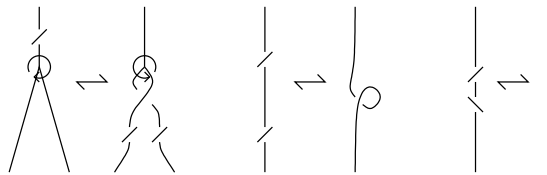
\begin{picture}(3264,1014)(5299,-2953)
\end{picture}%

\end{center}

The homomorphic expansion is constructed from $Z^{old}$ by first ignoring the twists, then placing
a factor of $e^{c/4}$
for each right twist and $e^{-c/4}$ for each left twist, where $c$ is a short (local) chord on the appropriate
strand (as well as the appropriate $\nu$ factors on dots and crosses).

\parpic[r]{\input{figs/nonorientedr.pstex_t}}
In this notation, the twisted theta graph representing $R$ is isotopic to the graph shown on the right. The
$Z$-value of the trivially embedded theta with three crosses is $1$, as shown before, so the value of this
graph is $e^{-c/4}$ on the first two strands and $e^{c/4}$ on the third, which is symmetric when 
pushed onto the first two strands: by the vertex invariance relation it equals 
$e^{{t_{12}}/2}$, as stated before.

\section{A note on the Kirby band-slide move and the LMO invariant}

In \cite{LMMO} Le, Murakami, Murakami and Ohtsuki construct an invariant
$\check{Z}$ of links which induces an invariant of $3$-manifolds, which
was recently {\it falsely} disputed (\cite{Ga}). The key 
step is proving that $\check{Z}$ is invariant under the Kirby
band-slide move $K2$, shown below:

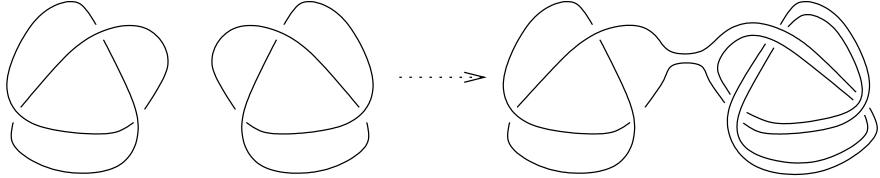
\begin{center}
 \begin{picture}(0,0)%
\includegraphics{figs/K2.pstex}%
\end{picture}%
\setlength{\unitlength}{4144sp}%
\begingroup\makeatletter\ifx\SetFigFont\undefined%
\gdef\SetFigFont#1#2#3#4#5{%
  \reset@font\fontsize{#1}{#2pt}%
  \fontfamily{#3}\fontseries{#4}\fontshape{#5}%
  \selectfont}%
\fi\endgroup%
\begin{picture}(5231,1059)(5126,-3580)
\end{picture}%

\end{center}

The problem is that this move is not a well-defined operation of links,
so somewhat cumbersome local considerations (``freezing'' local pictures 
or fixing bracketings) need to be used. $\check{Z}$ is defined to be a
normalized version of the classical Kontsevich integral $Z$, where an 
extra factor of $\nu$ is placed on each link component.

\parpic[r]{\input{figs/K2proof.pstex_t}}
Using our language, let us consider the sub-structure of dKTG' 
whose objects are links with possibly one $\theta$-graph component,
where circles are required to have two dots on them, and the 
$\theta$-graph component is required to have one dot on each edge. The
only operation we allow is unzipping edges of the $\theta$-graph.
We think of the theta as the two link components on which we want
to perform $K2$, ``fused together'' at the place where we perform
the operation. Unzipping the middle edge of the theta gives back the
original link (before $K2$), while unzipping a side edge produces
the link after $K2$ is performed. The invariance under $K2$ is then
a direct consequence of the homomorphicity of $Z$ with respect
to dotted unzip, as summarized by the figure on the right. Note that
in this case $Z$ is indeed $\check{Z}$ when restricted further to 
links, via replacing the dots by their values of $\nu^{1/2}$. In 
summary, we have proven the following:

\begin{theorem}\label{k2}
 There exists an expansion $\check{Z}$ for links with possibly one knotted
theta component, which is homomorphic with respect to unzipping any edge
of the theta component. When restricted to link (with no thetas), $\check{Z}$
agrees with the invariant $\check{Z}$ of \cite{LMMO}.
\end{theorem}

The reader may verify that the unzip property of Theorem \ref{k2} is
exactly the equivariance property required for the use of $\check{Z}$
in the construction of an invariant of rational homology spheres,
see \cite{LMO, BGRT1, BGRT2, BGRT3}.

\section{Appendix}

{\bf Proof of Theorem \ref{augm}.}
Let us denote the $n$-th piece of the classical finite type filtration by ${\mathcal F}_n$,
and the augmentation ideal by $\calI$. First we prove that $\calI={\mathcal F}_1$.

$\calI$ is linearly generated by differences, i.e., $\calI= \langle \gamma_1-\gamma_2 \rangle$, 
where $\gamma_1$ and $\gamma_2$ are KTGs of the same skeleton. ${\mathcal F}_1$ is linearly
generated by resolutions of $1$-singular KTGs, i.e. 
${\mathcal F}_1= \langle \gamma -\gamma' \rangle$, where $\gamma$ and $\gamma'$ differ in one
crossing change. Thus, it is obvious that ${\mathcal F}_1 \subseteq \calI$. The other direction,
$\calI \subseteq {\mathcal F}_1$ is true due to the fact that one can get to any knotting of a
given trivalent graph (skeleton) from any other through a series of crossing changes.

To prove that $\calI^n \subseteq {\mathcal F}_n$, we use that $\calI = {\mathcal F}_1$.
Recall that $({\mathcal F}_1)^n$ is generated by ``formulas'' containing $n$ $1$-singular KTGs, possibly
some further non-singular KTGs, joined by connected sums (the only binary operation), and possibly
with some other operations (unzips, deletes, orientation switches) applied. The connected sum
of a $k$-singular and an $l$-singular KTG is a $(k+l)$-singular KTG. It remains to check that 
orientation switch, delete and
unzip do not decrease the number of double points. Switching the orientation of an edge with a 
double point only introduces a negative sign.
Unzipping an edge with a double point on it produces a sum of two graphs with the same number
of double points. Deleting an edge with a double point on it produces zero. Thus, an element
in $({\mathcal F}_1)^n$ is $n$-singular, therefore contained in ${\mathcal F}_n$.

\parpic[r]{\begin{picture}(0,0)%
\includegraphics{figs/singular.pstex}%
\end{picture}%
\setlength{\unitlength}{3789sp}%
\begingroup\makeatletter\ifx\SetFigFont\undefined%
\gdef\SetFigFont#1#2#3#4#5{%
  \reset@font\fontsize{#1}{#2pt}%
  \fontfamily{#3}\fontseries{#4}\fontshape{#5}%
  \selectfont}%
\fi\endgroup%
\begin{picture}(1975,1590)(1902,-2223)
\end{picture}%
}
The last step is to show that ${\mathcal F}_n \subseteq \calI^n$, i.e., that one can write
any $n$-singular KTG as $n$ 1-singular, and possibly some further non-singular KTGs with a 
series of operations applied to them. The proof is in the same vein as proving that KTGs are
finitely generated \cite{Th}, as illustrated here by the example of a $2$-singular knotted theta-graph,
shown on the right. In the figures, a trivalent vertex denotes a vertex, while a $4$-valent
one is a double point. We start by taking a singular twisted tetrahedron for each double point,
a (non-singular) twisted tetrahedron for each crossing, and a standard tetrahedron for each
vertex, as shown in the figure below. We then apply a vertex connected sum (the composition 
of a connected sum and two unzips, as defined at the beginning of the proof of Theorem \ref{assoc})
along any tree connecting the tetrahedra, followed by sliding and unzipping edges, as shown below. 
The result is the desired KTG with an extra loop around it. Deleting the superfluous loop concludes the proof.
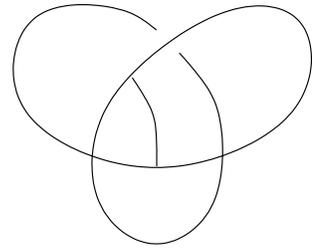
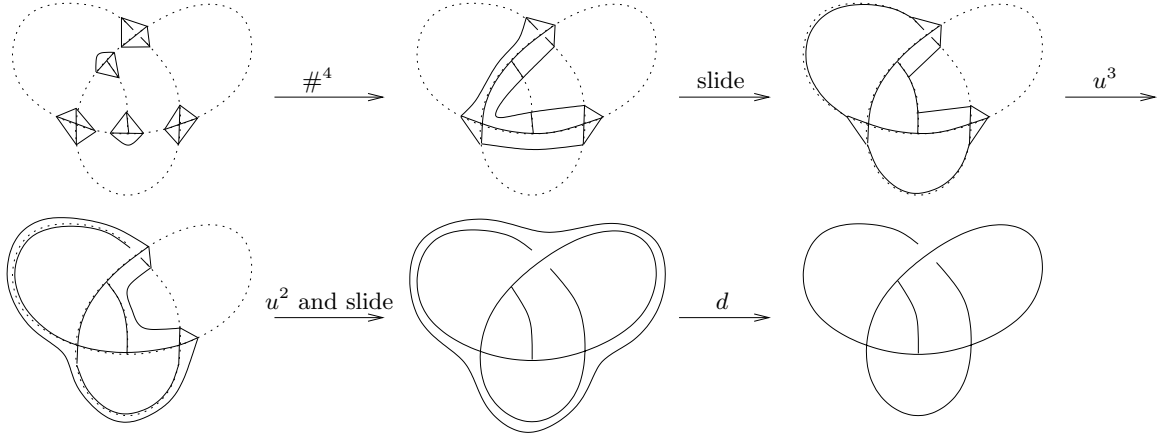
\begin{center}
 \input figs/construct.pstex_t
\end{center}
\vskip -8mm
\rightline{$\Box$}

\begin{remark}
\rm{Note that one is tempted to take connected sums along all dotted edges, but this would be a mistake:
connected sum can only connect separate connected components, adding an edge to a connected graph is
not well-defined; hence the slight awkwardness on choosing a ``spanning tree'' is necessary.}
\end{remark}

\bigskip

{\bf Proof of Theorem \ref{nonexistence}.}
Let us assume that a homomorphic universal finite type invariant of KTGs exists and
call it $Z$. By definition, $Z$ has to satisfy the following properties:

If $\gamma$ is a singular KTG (a KTG with finitely many transverse double
points), and $C_{\gamma}$ is its chord diagram (chords connect the pre-images
of the double points), then
\begin{center}
$Z(\gamma)=C_{\gamma}+$ higher order terms.
\end{center}

$Z$ commutes with all KTG operations, in particular, with orientation switch, unzip,
and connected sum: 
$$Z(S_e(\gamma))=S_e(Z(\gamma)),$$
$$Z(u_e(\gamma))=u_e(Z(\gamma)),$$
$$Z(\gamma \#_{e,f} \delta)=Z(\gamma)\#_{e,f}Z(\delta).$$

Let us denote the degree $k$ part of the values of $Z$ by $Z_k$.

To prove the theorem, we deduce a sequence of lemmas from the above properties, until we get a contradiction.

\begin{lemma}\label{nu}
Assume that $Z$ is homomorphic, as above, and let $Z(\circlearrowleft)=:\hat{\nu}$ be the value of the 
trivially framed unknot. Then,
in $\mathcal{A}(\circlearrowleft)$, $\hat{\nu}^2=\hat{\nu}$. This implies that all positive degree 
components of $\hat{\nu}$ are zero, i.e. $\hat{\nu}=1$.
\end{lemma}

\begin{proof}
$$Z(\circlearrowleft)=\raisebox{-5mm}{\begin{picture}(0,0)%
\includegraphics{figs/fig1.pstex}%
\end{picture}%
\setlength{\unitlength}{3947sp}%
\begingroup\makeatletter\ifx\SetFigFont\undefined%
\gdef\SetFigFont#1#2#3#4#5{%
  \reset@font\fontsize{#1}{#2pt}%
  \fontfamily{#3}\fontseries{#4}\fontshape{#5}%
  \selectfont}%
\fi\endgroup%
\begin{picture}(519,614)(1489,-818)
\put(1551,-566){\makebox(0,0)[lb]{\smash{{\SetFigFont{12}{14.4}{\rmdefault}{\mddefault}{\updefault}{\color[rgb]{0,0,0}$\hat{\nu}$}%
}}}}
\end{picture}%
}$$
Taking the connected sum of two unknots implies:
$$Z(\raisebox{-1mm}{\begin{picture}(0,0)%
\includegraphics{figs/fig2.pstex}%
\end{picture}%
\setlength{\unitlength}{3947sp}%
\begingroup\makeatletter\ifx\SetFigFont\undefined%
\gdef\SetFigFont#1#2#3#4#5{%
  \reset@font\fontsize{#1}{#2pt}%
  \fontfamily{#3}\fontseries{#4}\fontshape{#5}%
  \selectfont}%
\fi\endgroup%
\begin{picture}(596,230)(905,-296)
\end{picture}%
})=\raisebox{-5mm}{\input{figs/fig3.pstex_t}}$$
Unzipping the middle edge, we get an unknot back, which proves that $\hat{\nu}=\hat{\nu}^2$:
$$\raisebox{-5mm}{\begin{picture}(0,0)%
\includegraphics{figs/fig1.pstex}%
\end{picture}%
\setlength{\unitlength}{3947sp}%
\begingroup\makeatletter\ifx\SetFigFont\undefined%
\gdef\SetFigFont#1#2#3#4#5{%
  \reset@font\fontsize{#1}{#2pt}%
  \fontfamily{#3}\fontseries{#4}\fontshape{#5}%
  \selectfont}%
\fi\endgroup%
\begin{picture}(519,614)(1489,-818)
\put(1551,-566){\makebox(0,0)[lb]{\smash{{\SetFigFont{12}{14.4}{\rmdefault}{\mddefault}{\updefault}{\color[rgb]{0,0,0}$\hat{\nu}$}%
}}}}
\end{picture}%
}=Z(\circlearrowleft)
=\raisebox{-5mm}{\input{figs/fig5.pstex_t}} =\raisebox{-5mm}{\begin{picture}(0,0)%
\includegraphics{figs/fig4.pstex}%
\end{picture}%
\setlength{\unitlength}{3947sp}%
\begingroup\makeatletter\ifx\SetFigFont\undefined%
\gdef\SetFigFont#1#2#3#4#5{%
  \reset@font\fontsize{#1}{#2pt}%
  \fontfamily{#3}\fontseries{#4}\fontshape{#5}%
  \selectfont}%
\fi\endgroup%
\begin{picture}(519,614)(1489,-818)
\put(1531,-556){\makebox(0,0)[lb]{\smash{{\SetFigFont{10}{12.0}{\rmdefault}{\mddefault}{\updefault}{\color[rgb]{0,0,0}$\hat{\nu}^2$}%
}}}}
\end{picture}%
}$$
Since $\nu$ is an invertible element of $\mathcal{A}(\circlearrowleft)$, this implies
that $\nu=1$.
\end{proof}

Note that the conclusion of Lemma \ref{nu} implies that $\hat{\nu}=1$ as an element
of $\mathcal{A}(\uparrow)$, as $\mathcal{A}(\circlearrowleft) \cong \mathcal{A}(\uparrow)$.

\begin{corollary}\label{dumbbell}
$Z(\raisebox{-1mm}{\begin{picture}(0,0)%
\includegraphics{figs/dumbbell.pstex}%
\end{picture}%
\setlength{\unitlength}{4144sp}%
\begingroup\makeatletter\ifx\SetFigFont\undefined%
\gdef\SetFigFont#1#2#3#4#5{%
  \reset@font\fontsize{#1}{#2pt}%
  \fontfamily{#3}\fontseries{#4}\fontshape{#5}%
  \selectfont}%
\fi\endgroup%
\begin{picture}(654,194)(4984,-1958)
\end{picture}%
})=1 \in \calA(\uparrow_2)$.
\end{corollary}

\begin{proof}
The $Z$-value of the dumbbell graph can be viewed as an element of $\calA(\uparrow_2)$
by sweeping the middle edge free of chords, as we have done before. (In fact, this is not even
necessary, as it is an easy property of chord diagrams that any chord ending 
on a bridge makes a chord diagram zero.)

The statement of the corollary is obviously true, as the dumbbell graph is the connected sum of two unknots, 
and $Z$ commutes with connected sum.
\end{proof}

\begin{lemma}
 $\hat{\Phi}=Z\Bigg(\raisebox{-4mm}{\input{figs/tetr.pstex_t}}\Bigg)$ is an associator, 
with $\hat{R}=Z\Bigg(\raisebox{-4mm}{\input{figs/r2.pstex_t}}\Bigg)$.
\end{lemma}

\begin{proof}
 The proof is identical to the proof of Theorem \ref{assoc}, omitting all dots and crosses.
\end{proof}

\begin{corollary}\label{deg2}
$$Z_2\Bigg(\raisebox{-4mm}{\input{figs/tetr.pstex_t}}\Bigg)=
\alpha(t_{11}t_{23}-t_{33}t_{12})+\beta(t_{13}t_{12}-t_{13}t_{23})+\gamma(t_{12}t_{23}-t_{23}t_{12}),$$
where $t_{ij}$ denotes a single chord between strands $i$ and $j$ ($i$ and $j$ not necessarily different),
and $\alpha$, $\beta$ and $\gamma$ are complex numbers satisfying $\beta+\gamma=-\frac{1}{24}$. 
\end{corollary}

\begin{proof}
The universality of $Z$ (using that switching the crossing in $\hat{R}$ gives $(\hat{R}^{21})^{-1}$) implies
that the linear term
of $\hat{R}$ is $\frac{1}{2} t_{12}$.  

The statement then follows by solving the associator equations directly up to degree $2$, a 
straightforward computation that we omit here for the sake of concision.
\end{proof}

\begin{corollary}
 $Z_2(\raisebox{-1mm}{\begin{picture}(0,0)%
\includegraphics{figs/dumbbell.pstex}%
\end{picture}%
\setlength{\unitlength}{4144sp}%
\begingroup\makeatletter\ifx\SetFigFont\undefined%
\gdef\SetFigFont#1#2#3#4#5{%
  \reset@font\fontsize{#1}{#2pt}%
  \fontfamily{#3}\fontseries{#4}\fontshape{#5}%
  \selectfont}%
\fi\endgroup%
\begin{picture}(654,194)(4984,-1958)
\end{picture}%
})\neq 0$, in contradiction with Lemma \ref{dumbbell}.
\end{corollary}

\begin{proof}
 Switching the orientation of edge $1$ of the tetrahedron followed by unzipping the edge labeled $e$ below
results in a dumbbell graph:
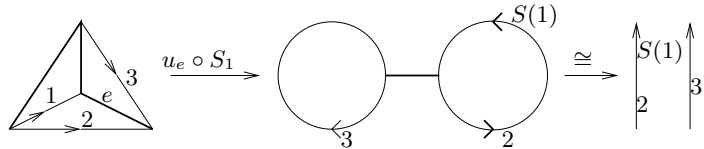
\begin{center}
 \input{figs/unziptetr.pstex_t}
\end{center}
Applying this to the degree $2$ term of Corollary \ref{deg2}, we get
$$-\frac{1}{24}t_{12}t_{11}+\alpha(t_{12}t_{11}+t_{11}t_{22})+\beta t_{12}^2+\gamma x,$$
where $x=\raisebox{-4mm}{\begin{picture}(0,0)%
\includegraphics{figs/chord.pstex}%
\end{picture}%
\setlength{\unitlength}{4144sp}%
\begingroup\makeatletter\ifx\SetFigFont\undefined%
\gdef\SetFigFont#1#2#3#4#5{%
  \reset@font\fontsize{#1}{#2pt}%
  \fontfamily{#3}\fontseries{#4}\fontshape{#5}%
  \selectfont}%
\fi\endgroup%
\begin{picture}(444,564)(6799,-1693)
\end{picture}%
}$. The reader can check that the result
is a nonzero element of $\calA(\uparrow_2)$, for all $\alpha$, $\beta$, $\gamma$ as above.
\end{proof}
This contradiction concludes the proof of Theorem \ref{nonexistence}. \qed

\end{document}